\documentclass{article}
\usepackage{amsmath,amstext}
\usepackage{amssymb}
\usepackage{amscd}
\usepackage{amsfonts}
\usepackage[english]{babel}
\usepackage[all]{xy}
\usepackage[latin1]{inputenc}
\numberwithin{equation}{section} \allowdisplaybreaks
\newtheorem{theorem}{Theorem}[section]
\newtheorem{lemma}[theorem]{Lemma}
\newtheorem{proposition}[theorem]{Proposition}
\newtheorem{corollary}[theorem]{Corollary}
\newtheorem{definition}[theorem]{Definition}

\newtheorem{remark}[theorem]{Remark}

\newcommand{\hess}{\mbox{\rm hess}}
\newcommand{\Pf}{\mbox{\rm Pf}}
\begin{document}

\bigskip

\begin{center}{\large Compatible structures on Lie algebroids\\ and
  Monge-Amp\`ere operators}\footnote{To be published in
  {\it{Acta. Appl. Math.}}, 2009.}
\end{center}

\bigskip

\begin{center} 
Yvette Kosmann-Schwarzbach and Vladimir Rubtsov
\end{center}

\bigskip

\begin{abstract}
We study pairs of structures, such as the Poisson-Nijenhuis structures,
on the tangent bundle of a manifold or, more generally, 
on a Lie algebroid or a Courant algebroid. 
These composite structures are defined by two of the following, a 
closed $2$-form, a Poisson bivector or a Nijenhuis tensor, with 
suitable compatibility assumptions. We establish the relationships
between $PN$-, $P \Omega$- and $\Omega N$-structures.
We then show that the non-degenerate Monge-Amp\`ere structures on
$2$-dimensional manifolds satisfying an integrability condition
provide numerous examples of such structures, 
while in the case of $3$-dimensional manifolds, such Monge-Amp\`ere
operators give rise to generalized complex structures or generalized
product structures on the cotangent bundle
of the manifold. 
\end{abstract}

\bigskip


\smallskip

Mathematics Subject Classification (MSC2000): 
53D17, 17B70, 58J60 (primary), 37K20, 37K25, 70G45 (secondary).

\bigskip

\bigskip

\begin{center}{\it For our friend Joseph Krasil'shchik on the occasion
    of his sixtieth birthday}
\end{center}

\section*{ Introduction}

On the tangent bundle of a manifold or, more generally, on
a Lie algebroid, we consider pairs of structures, 
such as the Poisson-Nijenhuis structures
which give rise to hierarchies
of Poisson structures (also called Hamiltonian structures) that play a
very important role
in the theory of completely integrable systems.
These structures are defined by 
closed $2$-forms, Poisson bivectors or $(1,1)$-tensors with vanishing Nijenhuis torsion.
When suitable compatibility assumptions are introduced, one obtains composite structures
called complementary $2$-forms, $PN$-, $P \Omega$- and $\Omega N$-structures.
Krasil'shchik contributed to the study of the algebraic nature of
Hamiltonian and bi-Hamiltonian structures and was the first to
underline the cohomological nature of their compatibility condition
(see \cite{K} and references therein). 
While the first part of this article is a comprehensive survey of the relationships 
between such composite structures and the related notion of Hitchin
pairs, the second part provides numerous examples 
arising from the theory of Monge-Amp\`ere equations. 

Our formulations and proofs make essential use of the 
{\it{big bracket}}, the even graded bracket on the space $\mathcal F$
of functions on the
cotangent bundle of a Lie algebroid considered as a supermanifold.
What we call the big bracket was first introduced by Kostant and Sternberg \cite{KS}; 
its use in the theory of Lie bialgebras is due to Lecomte and Roger 
\cite{LR} and was developed by one of us \cite{yks92}. Roytenberg 
extended it to Lie algebroids 
\cite{R} and Courant algebroids \cite{R2}. Recently, it has been used
by Antunes \cite{A} in the study of composite structures arising in
the theory
of sigma-models. In practice, all proofs are reduced to a
straightforward use of the graded Jacobi identity, sometimes 
repeatedly. While many of our results can be found in the literature
(see \cite{D}, which contains the references to earlier work 
by Magri, Gelfand and Dorfman,
Fokas and Fuchssteiner, see \cite{MM} \cite{MMR} \cite{KM},
and the more recent articles \cite{V1} \cite{V2} \cite{G} \cite{A}),
we claim
that our method unifies results, generalizing the known properties
from the case of manifolds 
to that of Lie algebroids and Courant algebroids. 
Our main argument is that the big bracket formalism
can be applied to problems in the geometric theory of partial
differential equations developed in \cite{Ly}\cite{LyR} \cite{KLR}
and \cite{B} \cite{BTh}. We also
stress that this theory can be considered in the general framework of
Lie algebroids, and we wish to introduce 
a general abstract theory of Monge-Amp\`ere structures on arbitrary
Lie algebroids.
In particular, the symplectic Monge-Amp\`ere equations defined by
$n$-forms on the cotangent bundle of a smooth, $n$-dimensional
manifold $M$ and, more generally, the Jacobi first-order systems,
defined by a set of $2$-forms on an $m+2$-dimensional manifold $M$,
can be viewed as ``deformations'' of the standard Lie algebroid
structure on the tangent bundle $T(T^*M)$ of $T^*M$. We shall indicate
some links between our approach and the approach to the geometric
structures
developed by Hitchin \cite{H} and Gualtieri \cite{Gua} in their studies
of generalized complex and K\"ahler structures, a new and fast
developing field of differential geometry.

In Section \ref{section1},
we introduce the big bracket, we recall the definition of Lie algebroids and  
give the explicit expression for the {\it{Dorfman bracket}}
on the double of a Lie bialgebroid
which is a {\it{derived bracket}} \cite{yks1996} \cite{yks2004}
of the big bracket. The Courant algebroid structure of the double of a
Lie bialgebroid is 
defined by the skew-symmetrized version of the Dorfman bracket, called
the {\it Courant bracket}. 
Section \ref{2} deals with general facts and formulas involving
bivectors, forms and $(1,1)$-tensors that will be used in subsequent
sections, 
and with Grabowski's formula \eqref{a} that expresses the Nijenhuis torsion of a
$(1,1)$-tensor in terms of the big bracket \cite{G}.
In Section \ref{representation}, we show that the adjoint actions of a
non-degenerate $2$-form and of its inverse bivector
induce a representation of ${\mathfrak{sl}}_2$ on $\mathcal F$, we
define the {\it primitive elements} and describe a Hodge-Lepage 
type decomposition  
of the elements in $\mathcal F$.

Section \ref{4} is a study of the {\it complementary $2$-forms} introduced
by Vaisman \cite{V1} \cite{V2}.
We prove that, given a Lie algebroid $A$, ``$\omega$ 
is a complementary $2$-form for 
the Poisson bivector $\pi$'' is a sufficient condition for the bracket
obtained by first dualizing the Lie algebroid structure of $A$
by $\pi$ and dualizing again by $\omega$ to be a Lie algebroid bracket
on $A$, whose expression we easily derive. 
A remark concerning the corresponding modular class (Section
\ref{modclass}) will be used in Section \ref{MA2modular}. 
Sections \ref{sectionPN} to \ref{table}
contain the detailed analysis of the 
structures introduced by Magri and Morosi \cite{MM} \cite{MMR} 
defined by a Poisson bivector and a Nijenhuis tensor, called
$PN$-structures (Section \ref{sectionPN}), by compatible
Poisson tensors (Section \ref{sectioncompatible}), by a closed
$2$-form and a Poisson bivector, called $P\Omega$-structures (Section
\ref{sectionpomega})
and by a closed $2$-form and a Nijenhuis tensor, called $\Omega
N$-structures (Section \ref{8}) and {\it Hitchin pairs} introduced by 
Crainic \cite{C} (Section \ref{9}).
A table and a diagram summarize the relationships between these
various structures.

Section \ref{courant} deals with Nijenhuis tensors on Courant
algebroid. We 
state Grabowski's theorem \cite{G} that characterizes {\it generalized
complex structures} 
by a simple equation in terms of the big bracket.

In Sections \ref{MA}-\ref{MA3}, we describe the geometry
of the symplectic {\it {Monge-Amp\`ere
equations}} and relate it to the structures discussed in the previous
sections, using the formalism of the big bracket.
Some of these results are reformulations of results in \cite{KLR} and
\cite{B} \cite{BTh}.
Section \ref{MA} introduces Monge-Amp\`ere structures on manifolds
and the associated Monge-Amp\`ere operators and equations. We recall
the definition of the {\it{effective forms}} 
and the one-to-one correspondence between
Monge-Amp\`ere operators and effective forms.
Section \ref{MA2} is devoted to the case of Monge-Amp\`ere structures
on $2$-dimensional manifolds, with an emphasis on the non-degenerate
case, when the Pfaffian of the defining $2$-form is nowhere vanishing.
We show that in the integrable case, {\it i.e.}, when the
Monge-Amp\`ere operator is equivalent to an operator with constant
coefficients,
the Monge-Amp\`ere structure gives rise to $PN$- and $\Omega
N$-structures and to a deformed Lie algebroid structure on $T(T^*M)$
which is unimodular.
More generally, a non-degenerate Monge-Amp\`ere structure of divergence type
defines a generalized almost complex structure on $T^*M$. If the
defining $2$-form is closed, this structure
is integrable and corresponds to a Hitchin pair.
The von Karman equation is an example where the integrability
condition is not satisfied and the associated composite structures do
not satisfy the compatibility condition. 
We then consider the first-order Jacobi differential systems which
generalize the Monge-Amp\`ere equations, and 
we describe the associated geometric structures on $2$-dimensional manifolds.
In Section \ref{MA3}, we proceed to study Monge-Amp\`ere operators on
$3$-dimensional manifolds, recall the classification of the
non-degenerate Monge-Amp\`ere operators, and we prove that when the
operator is non-degenerate, {\it i.e.}, when the Hitchin Pfaffian is
nowhere-vanishing, and has constant coefficients, there is
either an associated generalized complex structure or generalized
product structure on $T^*M$.
We conclude with a short discussion of two definitions of the
generalized Calabi-Yau manifolds.

\section{The big bracket}\label{section1}

When $A \to M$ is a vector bundle, let $T^*[2]A[1]$ denote the
cotangent bundle of the graded manifold $A[1]$ obtained from $A$ by
assigning degree $0$ to the coordinates on the base and degree $1$ to
the coordinates on the fibers.
The space $\mathcal F$ of smooth functions on $T^*[2]A[1]$
is a bigraded Poisson algebra \cite{R} \cite{yks2007}.
(See \cite{yks92} for the case where $M$ is a point and therefore $A$
is a vector space.)
If $(x^i, \xi^a)$, $i= 1, \ldots, {\rm dim}M$ and
$a = 1 , \ldots, {\rm rank} A$, are coordinates on $A[1]$, then
coordinates on $T^*[2]A[1]$ are
$(x^i, \xi^a, p_i, \theta_a)$, with bidegrees $(0,0), (0,1), (1,1),
(1,0)$, respectively. 
If an element $u$ 
of $\mathcal F$ is of bidegree
$(p+1,q+1)$, we call $|u| = p+q +2$ its (total) degree and 
we call $(p,q)$ its {\it{shifted bidegree}}, $p \geq -
1$, $q \geq - 1$. The space ${\mathcal F}^{p,q}$ of elements of
$\mathcal F$ of shifted bidegree $(p,q)$ contains the space of
sections of $\wedge^{p+1} A
\otimes \wedge^{q+1} A^*$. 

As the cotangent bundle of a graded manifold, $T^*[2]A[1]$ is
canonically equipped with an even Poisson structure. We denote the
even Poisson bracket on $\mathcal F$ by $\{~,~\}$, and we
call it the {\it big bracket}. The big bracket 
satisfies
$\{x^i,p_j\} = \delta^i_j$ and $\{\xi^a,\theta_b\} = \delta^a_b$,
so that $\{f,p_j\} = \partial_j f$, where $f \in C^\infty(M)$ and
$\partial_j f =\frac{\partial f}{\partial x^j}$.
This bracket is of bidegree
$(-1,-1)$ and of shifted bidegree $(0,0)$. It is skew-symmetric, 
$\{u,v\} = - (-1)^{|u|\,|v|} \{v,u\}$, for all $u$ and $v \in
{\mathcal F}$, 
and it satisfies the Jacobi
identity,
$$\{u,\{v,w\}\}= \{\{u,v\},w\} + (-1)^{|u| \, |v|} \{v, \{u,w\}\} \ ,
$$
for all $u$, $v$ and $w \in {\mathcal F}$. We often use the Jacobi
identity in the form,
$$\{\{u,v\},w\}= \{u,\{v,w\}\} + (-1)^{|v| \, |w|} \{\{u,w\}, v \} \ .
$$
The big bracket satisfies the Leibniz rule,
$$
\{u,v \wedge w \}= \{u,v\} \wedge w + (-1)^{|u| \, |v|} v \wedge  \{u,w\} \ ,
$$
or
$$
\{u \wedge v, w \}= u \wedge\{v, w\} + (-1)^{|v| \, |w|}  \{u,w\}
\wedge v \ .
$$

The space of sections of a vector bundle $E$ is denoted by 
${\rm{\Gamma}} E$. We call a section of $\wedge^\bullet E$ (resp.,
$\wedge^\bullet E^*$) a multivector (resp., a form) on $E$.
Accordingly, we use the terms vector, bivector, $k$-form,
$(p,q)$-tensor, 
etc. All manifolds and maps are assumed to be smooth.

\subsection{Lie algebroids}
A {\it Lie algebroid} structure on $A \to M$ is an element $\mu$ of
$\mathcal F$ of shifted bidegree $(0,1)$ such that
$$
\{\mu,\mu\} = 0 \ .
$$
The {\it Schouten bracket} of
multivectors, {\it i.e.}, sections of $\wedge^\bullet A$, $X$ and $Y$, is
$$
[X,Y]_\mu =  \{\{X,\mu\},Y\} \ .
$$
In particular, this formula defines
the Lie bracket of $X$ and $Y \in {\rm{\Gamma}} A$ as well as
the anchor of $A$, $\rho: A \to TM$, by
$$
\rho(X)f = \{\{X,\mu\},f\} \ ,
$$
for $X \in {\rm{\Gamma}} A$ and $f \in C^{\infty}(M)$.

The {\it Lie algebroid differential}
acting on sections of $\wedge^\bullet A^*$ is denoted by $d_\mu$, thus
$$
d_\mu =
\{\mu, \cdot \} \ .
$$
The Lie derivative of forms by $X \in {\rm{\Gamma}}
A$ is defined to be the graded commutator, 
${\mathcal L}^\mu_X = [i_X, d_\mu]$.

A {\it Lie bialgebroid} is defined by $\mu \in \mathcal F^{0,1}$
and $\gamma \in \mathcal F^{1,0}$ such that $\{\mu +
\gamma , \mu + \gamma \} = 0$. More generally, a {\it proto-bialgebroid}
is defined by $S = \phi + \mu+ \gamma + \psi$, where $\psi \in
{\rm{\Gamma}} (\wedge^3 A^*)$
and $\phi \in {\rm{\Gamma}} (\wedge^3 A)$, such that $\{S,S\}=0$.
{\it Lie quasi-bialgebroids} correspond to $\psi =0$, while 
{\it quasi-Lie bialgebroids} correspond to $\phi = 0$.

\subsection{The Dorfman bracket}\label{sectiondorfman}

If $(A, \mu,\gamma)$ is a Lie bialgebroid,  its double is the vector
bundle, $A \oplus A^*$, equipped with
the {\it Dorfman bracket} defined by
\begin{equation}\label{dorfman}
[u,v]_D =  \{\{u,\mu + \gamma \},v\} \ ,
\end{equation}
for $u$ and $v \in {\rm{\Gamma}}(A \oplus A^*)$.
The skew-symmetrized Dorfman bracket is called the {\it {Courant bracket}} 
and $A \oplus A^*$ with the Dorfman bracket is a {\it Courant
algebroid}. Since the Dorfman bracket (see \cite{D} \cite{yks2004}
\cite{G}) 
is a derived bracket, it is a Loday-Leibniz bracket and therefore satisfies the
(graded) Jacobi identity in the sense that, for each $u \in
{\rm{\Gamma}}(A\oplus A^*)$, $[u, \cdot]_D$ is a derivation of 
the bracket $[~, ~]_D$ (see \cite{yks1996} \cite{yks2004}).
More generally, Formula \eqref{dorfman} defines a Loday-Gerstenhaber
bracket on ${\rm{\Gamma}}(\wedge^\bullet A \otimes \wedge^\bullet A^*)$.

Explicitly,
for $X \in {\rm{\Gamma}} A$ and $\alpha \in {\rm{\Gamma}}(A^*)$,
$$
[X,\alpha]_D =  \{\{X,\mu\},\alpha\} + \{\{X,\gamma\},\alpha\}
= \{X, \{\mu, \alpha\} \}+ \{\mu,
\{X,\alpha\}\} - \{\{\gamma,X\},\alpha\}
$$
$$
= i_X(d_\mu\alpha) + d_\mu(i_X\alpha) - i_\alpha (d_\gamma X)
= {\mathcal L}^\mu_X\alpha - i_\alpha (d_\gamma X)  \ ,
$$
while
$$
[\alpha,X]_D = \{\{\alpha,\mu\}, X\} + \{\{\alpha,\gamma\}, X\}
= - \{\{\mu, \alpha\}, X\} + \{\alpha, \{\gamma,X\}\} + \{\gamma, \{\alpha,X\}\}
$$
$$
= - i_X(d_\mu\alpha) +i_\alpha(d_\gamma X) + d_\gamma(i_\alpha X)
= {\mathcal L}^\gamma_\alpha X  - i_X(d_\mu\alpha) \ .
$$
Therefore, for $X$ and $Y \in {\rm{\Gamma}} A$, $\alpha$ and $\beta \in
{\rm{\Gamma}}(A^*)$,
\begin{equation}
[X+\alpha, Y+\beta]_D = [X,Y]_\mu + {\mathcal L}^\gamma_\alpha Y
- i_\beta (d_\gamma X)
+ [\alpha,\beta]_\gamma
+ {\mathcal L}^\mu_X\beta - i_Y(d_\mu\alpha)  \ .
\end{equation}

In the case of the
standard Courant algebroid, $TM \oplus T^*M$, by assumption, $\gamma =0$ and $d_\mu$
is the de Rham differential, $d$. Thus,
for $X \in {\rm{\Gamma}}(TM)$ and $\alpha \in {\rm{\Gamma}}(T^*M)$,
$$
[X,\alpha]_D = \{\{X,\mu\},\alpha\} = \{X, \{\mu, \alpha\}\} + \{\mu,
\{X,\alpha\}\}
= i_X(d\alpha) + d(i_X\alpha) = {\mathcal L}_X\alpha \ ,
$$
and
$$
[\alpha,X]_D = \{\{\alpha,\mu\}, X\} = - \{\{\mu, \alpha\}, X\}
= - i_X(d\alpha) \ .
$$
In addition, it is clear that, for vector fields $X$ and $Y$, $[X,Y]_D$
is the Lie bracket, and for $1$-forms, $\alpha$ and $\beta$,
$[\alpha,\beta]_D = 0$, this bracket vanishes on pairs of $1$-forms. Therefore
\begin{equation}
[X+\alpha, Y +\beta]_D = [X,Y] + {\mathcal L}_X\beta-i_Y(d\alpha) \ .
\end{equation}
We compute these brackets on $T^*[2]TM[1]$ in local coordinates,
$(x^i,\xi^i,p_i,\theta_i)$.
Here $\mu = p_i\xi^i$. Let $X =X^i\theta_i$ and $\alpha =
\alpha_i\xi^i$. Then
$$[X,\alpha]_D = \{\{X^i\theta_i, p_j\xi^j\}, \alpha_k\xi^k\}
= \{X^i p_i - \partial_j X^i \theta_i \xi^j, \alpha_k \xi^k\}
= X^i \partial_i \alpha_k \xi^k + \partial_j X^i \alpha _i \xi^j \ ,
$$
which is the expression of ${\mathcal L}_X\alpha$ in local coordinates.
Similarly,
$$
[\alpha,X]_D =\{\{ \alpha_i\xi^i, p_j\xi^j \},  X^k\theta_k\}
= \{-\partial_j \alpha_i \xi^i \xi^j , X^k \theta_k\}
= - X^k \partial_k \alpha_i \xi^i + X^k \partial_j \alpha_k \xi^j \ ,
$$
which is the expression of $-i_X(d\alpha)$ in local coordinates.

\begin{remark}{\rm{When $\mu$ is replaced by $\mu + H$, where $H$ is
    a $d_\mu$-closed $3$-form, the equation $\{\mu+H,\mu+H\} =0$ is
    satisfied, and one obtains the {\it Dorfman bracket with
      background}, $[~,~]_{D,H}$, on
    ${\rm{\Gamma}}(A\oplus A^*)$,
\begin{equation}
[X+\alpha, Y+\beta]_{D,H} = [X+\alpha, Y+\beta]_D + i_{X\wedge Y}H \ ,
\end{equation}
making $A\oplus A^*$ a {\it twisted Courant algebroid} \cite{SW} \cite{R} \cite{Gua} .

\medskip

Any $2$-form $B$ on $A$ defines a gauge transformation,
$\widehat B: X+\alpha \mapsto X+\alpha +i_XB$, satisfying
$$
[\widehat B(X+\alpha), \widehat B(Y+\beta)]_{D,H}=
\widehat B([X+\alpha, Y +\beta]_{D,H-d_\mu B}) \ .
$$
If $B$ is $d_\mu$-closed, then $\widehat B$ is an automorphism of $(A
\oplus A^*, [~,~]_{D,H})$.}}
\end{remark}

\section{Tensors and the big bracket}\label{2}
We shall need various preliminary results
concerning tensors on a Lie algebroid.

\subsection{Bivectors, forms and $(1,1)$-tensors}
Let $\pi^{\sharp} : A^{*} \to A$ be the map defined by a bivector
$\pi$, where $\pi^\sharp
\alpha = i_\alpha \pi$, for $\alpha \in {\rm{\Gamma}}(A^*)$. Then
\begin{equation}\label{piflat}
\pi^\sharp
\alpha = \{\alpha, \pi\} \ .
\end{equation}

Let $\omega^{\flat}:
A\to A^{*}$ be the map defined by a $2$-form $\omega$, where
$\omega^\flat X = - i_X \omega$, for $X \in {\rm{\Gamma}} A$.
Then
\begin{equation}
\omega^\flat X = \{\omega, X \} \ .
\end{equation}

Let $N\,^{\hat{}}
 : {\rm{\Gamma}} A \to {\rm{\Gamma}} A$ be the linear map induced by a vector bundle
endomorphism of $A$. Then $N\,^{\hat{}}$
can be identified with a $(1,1)$-tensor on
$A$, more precisely with a section $N$ of $A^*\otimes A$, by setting
\begin{equation}\label{defN}
N\,^{\hat{}}(X)=\{X,N\} \ ,
\end{equation}
for all $X\in {\rm{\Gamma}} A$.
In local coordinates, if $N\,^{\hat{}}$ has
  components $N^a_b$, then $N = N^a_b \xi^b \theta_a$.
We shall not distinguish between $N\,^{\hat{}}$ and $N$,
and we shall abbreviate $N\,^{\hat{}}$ to $N$.

\begin{lemma}\label{lemma2} The map
$N= \pi^{\sharp} \circ \omega^{\flat} : A\to A$ considered as a
section of $A^* \otimes A$ is
$$
N = \{\pi, \omega\} \ .
$$
\end{lemma}

\noindent{\it Proof}
By the Jacobi identity, since $\{ X,\pi\} = 0$, for all $X \in {\rm \Gamma}A$,
$$
\{X,N\}= \{X, \{\pi, \omega\}\} = \{\{ \omega, \pi \},X \} =\{\{
\omega, X \}, \pi \} =
 \{\omega^\flat X, \pi\}
= (\pi^\sharp \circ \omega^\flat) (X) \ .
$$
This proves the result, in view of \eqref{defN}. \hfill $\square$

\medskip

In local coordinates, let
 $\pi = \frac{1}{2} \pi^{ab} \theta_{a} \theta_{b}$ and
$\omega = \frac{1}{2} \omega_{ab} \xi^{a}\xi^{b}$.
For $\alpha = \alpha_a \xi^a \in {\rm{\Gamma}}(A^*)$, $
\pi^\sharp \alpha = \pi^{ab} \alpha_a
\theta_b$. For $X = X^a \theta_a \in {\rm{\Gamma}} A$,
$\omega^\flat X = \omega_{ab}X^b \xi ^a$, whence
$(\pi^{\sharp}\circ \omega^{\flat}) (X) = \pi^{ab} \omega_{ac}
X^{c} \theta_{b}$.
On
the other hand,
$$
\{ \pi, \omega \}=
\frac{1}{4} \pi^{ab} \omega_{cd} \{ \theta_{a}\theta_{b} ,
\xi^{c} \xi^{d}\} = \pi^{ab} \omega_{bc} \theta_{a} \xi^{c} \ ,
$$
whence $\{X , \{\pi, \omega\}\} =
\pi^{ab} \omega_{ac}
X^{c} \theta_{b}$.

\medskip

In particular,  if $\pi$ is non-degenerate, and if $\pi$ and $\omega$
are inverses of one another, by
definition, $\pi^\sharp \circ \omega^\flat = {\mathrm{Id}}_A$
and $\omega^\flat (X)= - i_X\omega$, for all $X \in {\rm{\Gamma}} A$.
Then the $(1,1)$-tensor
$\{\pi, \omega\}$ is the identity of $A$, ${\mathrm{Id}}_A$.
In Section \ref{representation} below, 
we denote the adjoint action of ${\mathrm{Id}}_A$,
$\{{\mathrm{Id}}_A, \cdot\}$, by ${\bf I}$.

In local coordinates, $\pi^{ab}\omega_{ac} = \delta^b_c$ and
$\{\pi, \omega\} = \xi^a\theta_a$, satisfying $\{X,
\xi^a\theta_a\}= X$, for all $X \in {\rm{\Gamma}} A$.

This relation is a particular case of a general result, proved in
\cite{yks2007}: when  $\pi$ and $\omega$ are inverses
of one another, for $u \in \mathcal F^{p,q}$,
in particular for $u\in {\rm{\Gamma}} (\wedge^{p+1}A
\otimes \wedge^{q+1}A^*)$,
\begin{equation}\label{inverse}
\{\{\pi, \omega\}, u\} = \{{\mathrm{Id}}_A , u\} = (q-p) u \ ,
\end{equation}
or, in local coordinates,
\begin{equation}\label{inverse2}
\{ \xi^a\theta_a , u \} = (q-p)u \ .
\end{equation}
Therefore, for a bracket $\mu$ (resp., cobracket $\gamma$) on $A$,
$$
\{ {\rm {Id}}_A, \mu\} = \mu \quad \quad \quad {\mathrm{(resp.,}}
\quad 
\{ {\rm {Id}}_A, \gamma\} = -
\gamma)
$$ and for a $3$-form  $\psi$ (resp., $3$-tensor $\phi$) on $A$,
$\{ {\rm {Id}}_A, \psi\} = 3 \psi$ (resp., $\{ {\rm {Id}}_A, \phi\} = -
3 \phi$).

\subsection{Deformed brackets and torsion}

Let $(A,\mu)$ be a Lie algebroid.
Let $N \in {\rm{\Gamma}}(A^* \otimes A)$ be a $(1,1)$-tensor on $A$, an element
of shifted bidegree $(0,0)$.
Then the {\it{deformed structure}},
$$
\mu_N = \{N,\mu\} \ ,
$$
defines an anchor $\rho \circ N$ and a
skew-symmetric bracket on $A$ which we shall denote by $[~,~]^\mu_N$,
Explicitly,
\begin{equation}\label{deformedbracket}
[X,Y]^\mu_N = \{\{X,\{N,\mu\}\} ,Y\} \ ,
\end{equation}
for $X$ and $Y \in {\rm{\Gamma}} A$.

\begin{lemma}\label{bracketN}
The bracket
$[~,~]^{\mu}_{N}$ is such that,
for $X$, $Y \in
{\rm{\Gamma}} A$,
\begin{equation}\label{deformed}
[X,Y]^\mu_{N} = [NX,Y]_{\mu} + [X,NY]_{\mu} - N
[X,Y]_{\mu} \ .
\end{equation}
\end{lemma}

\noindent{\it Proof} By definition,
$$
[X,Y]^\mu_{N}= \{\{X,\{N,\mu\}\},Y\} = \{\{\{X,N \},\mu\},Y\} +
\{\{N,\{X,\mu\}\},Y\}
$$
$$
= [NX,Y]_\mu + \{N,\{\{X,\mu\},Y\}\} + \{\{N,Y\},\{X,\mu\}\}
$$
$$
= [NX,Y]_\mu + [X,NY]_\mu - N[X,Y]_\mu \ ,
$$
where we have used the Jacobi identity and the definition of $[~,~]_\mu$.
\hfill $\square$

\medskip

The bracket $[~,~]^\mu_N$ is called the {\it deformed}
\cite{KM} (or
{\it contracted} \cite{CGM1} \cite{CGM2} \cite{G}) bracket of $[~,~]_\mu$.
By ${\mathcal T}_\mu N$ we denote the {\it Nijenhuis torsion} of $N$
defined by
\begin{equation}
({\mathcal T}_\mu N)(X,Y) =
[NX,NY]_\mu - N ([NX,Y]_\mu + [X,NY]_\mu) + N^2 [X,Y]_\mu \ ,
\end{equation}
for all $X$ and $Y \in {\rm{\Gamma}} A$. It is clear that
$({\mathcal T}_\mu N)(X,Y) =
[NX,NY]_\mu - N ([X,Y]^\mu_N)$.

\begin{proposition}\label{Nijtorsion}
In terms of the big bracket,
\begin{equation}\label{a}
{\mathcal T}_\mu N = \frac{1}{2} (\{N, \{N, \mu\}\} - \{N^2 ,\mu\}) \ ,
\end{equation}
and
\begin{equation}\label{b}
\frac{1}{2} \{\{N,\mu\},\{N,\mu\}\} = \{\mu,{\mathcal T}_\mu N \} \ .
\end{equation}
\end{proposition}

\noindent{\it{Proof}} See Grabowski \cite{G} for Formula
\eqref{a}, which is proved by a simple calculation.
Formula \eqref{b} follows from \eqref{a} by an application of the
Jacobi identity. \hfill $\square$

\medskip

\begin{remark}{\rm{Formula \eqref{a} 
can also be viewed as a particular case of
Formulas (5.22) and (5.16) of \cite{yks1996}, taking into account the
fact that, for vector-valued forms, the
big bracket and the Richardson-Nijenhuis bracket coincide up to sign
(see \cite{yks92}), or as a particular case of Formula (3.14) of \cite{yks2004}.
Formula \eqref{a} also appears in a slightly different form
in \cite{D}, Section 3.3. It plays an essential role in \cite{A}.}}
\end{remark}
\subsection{Nijenhuis structures}\label{sectionNij}

Let $N \in {\rm{\Gamma}}(A^* \otimes A)$ be a $(1,1)$-tensor on $A$, thus $N$
is an element
of shifted bidegree $(0,0)$.
Then the {\it{deformed structure}} bracket is defined by
\eqref{deformedbracket}, and its explicit expression is 
Formula \eqref{deformed} above.
We have denoted the {\it Nijenhuis torsion} of $N$ by ${\mathcal T}_\mu N$.
The following result (see, {\it e.g.}, \cite{KM}) is an immediate corollary
of Proposition \ref{Nijtorsion}.

\begin{theorem}\label{CNS}
A necessary and sufficient condition (resp., a sufficient condition) for
the deformed structure $\mu_N = \{N,\mu\}$ to be a Lie
algebroid structure on $A$ is
$$
\{\mu, {\mathcal T}_\mu N \} = 0 \ .
$$
(resp., ${\mathcal T}_\mu N = 0$).
\end{theorem}

When $\{\mu,{\mathcal T}_\mu N\}=0 $, we call
${\mathcal T}_\mu N$ a $d_\mu$-{\it {cocycle}}.

\begin{remark}\label{remarkcompatible}{\rm{The deformed Lie algebroid structure
$\mu_N$ is  compatible with $\mu$
in the sense that $\mu + \mu_N$ is a Lie
algebroid structure, {\it i.e.},
$\{\mu + \mu_N, \mu +\mu_N\} =0$.}}
\end{remark}

 The operator on ${\rm{\Gamma}}(\wedge^\bullet A^*)$ associated to $\mu_N$
is $d_{\mu_N} = \{\mu_N, \cdot \} = \{\{N,\mu\}, \cdot \}$.

\medskip

By definition, a $(1,1)$-tensor $N$ is an {\it almost
complex structure} if $N^2= -
{\mathrm{Id}}_A$, and an almost complex structure
$N$ is a {\it complex structure} if ${\mathcal T}_\mu N = 0$.

\begin{proposition}
An almost complex structure $N$ is a complex structure if and only if
\begin{equation}\label{complex}
\{\{N,\mu\},N\} = \mu \ .
\end{equation}
\end{proposition}

\noindent{\it Proof}
Equation \eqref{complex}
follows from \eqref{a} and the relation $\{ {\rm Id_A}, \mu \} = \mu$, a
particular case of \eqref{inverse}.
\hfill $\square$

\subsection{Bivectors and $3$-forms}\label{bivectors}
\begin{lemma}\label{lemma3}
If $\pi$ is a bivector and $\psi$ is a $3$-form on $A$,
then, for $X$, $Y \in{\rm{\Gamma}} A$,
\begin{equation}\label{3form}
\{\{ X , \{ \psi, \pi \}\},Y\}
= \pi^{\sharp} (i_{X\wedge Y} \psi) \ ,
\end{equation}
with the convention $i_{X\wedge Y}=i_X \circ i_Y$.
\end{lemma}

\noindent{\it Proof} Since $\{ X,\pi\} = \{ Y,\pi\}
= 0$, by the Jacobi identity, 
$$
\{\{ X , \{ \psi, \pi \}\} ,Y\} = \{\{ \{ X , \psi \} ,\pi\} ,Y\} =
\{\{\{ X,\psi \},Y\},\pi\}\ .
$$
Now $\{ X,\psi\} = i_{X} \psi$, $\{\{ X,\psi \} ,Y\} = i_{X}
i_{Y}\psi$. Applying \eqref{piflat} to the $1$-form $\alpha
=\{\{X,\psi\},Y\}$, we obtain \eqref{3form}.
\hfill $\square$

\medskip

In local coordinates, let $\pi = \frac{1}{2} \pi^{ab} \theta_{a}\theta_{b}$ and
$\psi = \frac{1}{6} \psi_{abc} \xi^{a}\xi^{b} \xi^{c}$. Then $\{
\psi, \pi\} = \frac{1}{2} \pi^{cd} \psi_{abc}  \xi^{a}
\xi^{b} \theta_d$. For $X = X^{a}\theta_{a}$, $\{ X , \{ \psi,\pi\}\}
= \pi^{cd}
\psi_{abc} X^{a} \xi^{b} \theta_d$, and for $Y = Y^{a}\theta_{a}$,
$\{\{ X , \{ \psi ,\pi\}\} ,Y\} = \pi^{dc} \psi_{abc} X^{a} Y^{b}
\theta_{d}$. On the other hand,
$
i_{X} i_{Y} \psi = \psi_{abc} X^{b} Y^{a} \xi^{c}
$,
and $\pi^{\sharp} (i_{X \wedge Y} \psi) = \pi^{dc} \psi_{abc} X^a
Y^b \theta_{d}$.

\section{A representation of ${\mathfrak{sl}}_2$}\label{representation}

For $u\in{\mathcal F}^{p,q}$, we shall 
call $w(u)=q-p$ the {\it weight} of $u$.
Let $\pi$ be a non-degenerate bivector and $\omega$ a $2$-form
on $A$ which are inverses of one another. Then
$\{\pi,\omega\}={\rm{Id}}_A$. Set
${\mathrm{ad}}_\omega = \{\omega, \cdot \}$ and
${\mathrm {ad}}_\pi = \{\pi, \cdot \}$. Then $ {\bf I} =
\{\{\pi, \omega, \}, \cdot \}$ acts on ${\mathcal{F}}$ by
$$
{\bf{I}} (u) = w(u)u \ ,
$$
for $u \in {\mathcal{F}}^{p,q}$ (see Formula \eqref{inverse}).
Let ${\mathrm{ad}}'_\pi = \{.,\pi\} = - {\mathrm{ad}}_\pi$
be the right adjoint action of $\pi$. Then
$$
[{\bf I} , {\mathrm{ad}}_\omega] = 2 {\mathrm{ad}}_\omega \ ,  \quad  \quad   [{\bf I} ,
{\mathrm{ad}}'_\pi] = - 2 {\mathrm{ad}}'_\pi \
\quad  \quad  [ {\mathrm{ad}}_\omega, {\mathrm{ad}}'_\pi]
= {\bf I} \ ,
$$
where $[~,~]$ denotes the commutator of operators.
Therefore the operators $({\mathrm{ad}}_\omega , {\mathrm{ad}}'_\pi , \bf I)$
define a representation of
$\mathfrak{sl}_2$ on the linear space ${\mathcal{F}}$ which restricts to the linear space of all tensors, analogous to the
representation on forms in \cite{KLR}. Then
$$
{\mathrm{ad}}_{\omega}
({\mathcal{F}}^{p,q}) \subset {\mathcal{F}}^{p-1,q+1},\ {\mathrm{ad}}'_{\pi}({\mathcal{F}}^{p,q})\subset {\mathcal{F}}^{p+1,q-1} \ ,
$$
so that $w({\mathrm{ad}}_{\omega}u) = w(u)+2$ and $w({\mathrm{ad}}'_{\pi}u)= w(u)-2$.

\begin{definition}\label{prim}
An element $u\in {\mathcal{F}}^{p,q}$ is called {\emph {primitive}}
if $u\in \ker ({\mathrm{ad}}'_\pi)$, {\it i.e.}, if $\{u,\pi\} = 0.$
\end{definition}

The next statement follows from the definitions.
\begin{lemma}\label{intersect}
For $\pi$ and $\omega$ inverses of one another,
$\ker({\mathrm{ad}}'_\pi)\cap \ker({\mathrm{ad}}_\omega) \subset
\oplus_{p \geq - 1}{\mathcal{F}}^{p,p}$.
\end{lemma}

The inverse inclusion is not valid since counter-examples are
furnished by $(1,1)$-tensors $N$ of
shifted bidegree $(0,0)$ such that the $2$-form
$\{{\mathrm{ad}}_{\omega},N\}$ or the bivector
$\{{\mathrm{ad}}'_{\pi},N\}$ does not vanish, {\it e.g.}, when $N$ 
is a multiple of the identity.

The following theorem is an analogue of the Hodge-Lepage decompositions in K\"ahler \cite{W} and symplectic
\cite{LM} \cite{KLR} geometry. We first prove a lemma.

\begin{lemma}\label{weight}
Let $u\in {\mathcal{F}}$ be of weight $w(u)$. Then, for any $k\geq 0$,
$$
{\bf I}({\mathrm{ad}}_\omega^ku) = (w(u) +2k){\mathrm{ad}}_\omega^k u \ .
$$
If u is primitive, then
$$
 \quad 
\quad  \quad  \quad  \quad  \quad  \quad  \quad  \quad
{\mathrm{ad}}'_\pi ({\mathrm{ad}}_\omega ^k u) = - k(w(u) + k -1)
{\mathrm{ad}}_\omega^{k-1} u \ . \quad  \quad  \quad  \quad  \quad
\quad  \quad  
\quad \quad 
\hfill \square
$$
\end{lemma}

\noindent{\it {Proof}}
The first formula follows from $w({\mathrm{ad}}_\omega ^k u)=w(u)+2k$.
The second is proved by recursion on $k$.
\hfill $\square$

\medskip

From the complete reducibility of finite-dimensional representations
of semi-simple Lie algebras, and from
Lemma \ref{weight} we obtain the
following result \cite{Ly} \cite{KLR}.
\begin{theorem}\label{hodge}
Any element $u\in {\mathcal{F}}^{p,q}$ admits the decomposition,
\begin{equation}\label{hll}
u = u_0 + {\mathrm{ad}}_\omega u_1 + {\mathrm{ad}}_\omega^2 u_2
+\ldots +{\mathrm{ad}}_\omega^k u_k +\ldots \ ,
\end{equation}
where each $u_k, k\geq 0$, is a uniquely defined primitive element of
${\mathcal{F}}^{p+k,q-k}$ of weight $w(u)- 2k$.
\end{theorem}

\section{Complementary $2$-forms for Poisson structures}\label{4}

The complementary $2$-forms with respect to a Poisson
structure on a Lie algebroid were defined and studied
by Vaisman \cite{V1} \cite{V2}. We shall describe the complementary
$2$-forms on a Lie algebroid $A$ and their properties 
by means of the big bracket on $\mathcal F$.
The method of proof using the big bracket gives a clear view of their
nature and properties.

\subsection{Poisson bivectors}\label{section4.1}

We recall several well known facts 
concerning the Poisson structures on Lie algebroids \cite{R} \cite{PM232}.
\begin{lemma}
Let $(A,\mu)$ be a Lie algebroid.
If $\pi \in {\rm{\Gamma}} (\wedge^{2}A)$, then
$${\gamma_\pi  = \{
  \pi,\mu\}}$$
is of shifted bidegree $(1,0)$, and $\gamma_\pi $ is a Lie
algebroid structure on  $A^{*}$ if and only if
\begin{equation}
\{ \gamma_\pi  ,\gamma_\pi \}=0 \ .
\end{equation}
\end{lemma}

The next lemma gives conditions for the construction, 
from a bivector on a Lie algebroid, of a Lie
algebroid structure on the dual vector bundle.
Since $\{\mu ,\mu\} =0$,
$$\{ \gamma_\pi ,
\gamma_\pi \} = \{\{ \pi,\mu\} , \{ \pi,\mu\}\} = \{\{\{
\pi,\mu\} ,\pi\} ,\mu\} = \{ [\pi,\pi]_{\mu},\mu\} \ .
$$
Therefore
\begin{lemma} A necessary and sufficient condition for $\gamma_\pi $ to be
a Lie algebroid structure on $A^{*}$ is
\begin{equation}\label{1}
\{ \mu, [\pi,\pi]_{\mu} \} = 0 \ ,
\end{equation}
while a sufficient condition is
\begin{equation}
[\pi,\pi]_{\mu}=0 \ ,
\end{equation}
{\it i.e.}, $\pi$ is a {\it Poisson bivector}.
\end{lemma}

The bracket defined by $\gamma_\pi = \{\pi,\mu\}$ on
${\rm{\Gamma}}(\wedge^\bullet A^*)$
is usually denoted simply by $[~,~]_\pi$. Thus, by definition, 
$$
\{\{ \alpha , \{
\pi,\mu\}\},\beta\} = [\alpha ,\beta]_{\pi} \ ,
$$
for all $\alpha$ and $\beta \in {\rm {\Gamma}}( \wedge^\bullet A^*)$.
The following lemma is proved in \cite{R} \cite{PM232}.

\begin{lemma}\label{lemma1}
The bracket defined by $\gamma_\pi = \{\pi,\mu\}$ on
${\rm{\Gamma}}(\wedge ^\bullet A^*)$
is the {\it Koszul bracket} of forms. In
particular, for all $f\in C^{\infty} (M)$,
$\alpha ,\beta \in {\rm{\Gamma}} ( A^{*})$,
$$
\{\{ \alpha , \{ \pi,\mu\}\}, f\} =
((\rho \circ \pi^\sharp) \alpha) \cdot f\ ,
$$
$$
\{\{ \alpha , \{ \pi,\mu\}\}, \beta\} = {\mathcal L}^\mu_{\pi^\sharp
  \alpha} \beta - {\mathcal L}^\mu_{\pi^\sharp
  \beta} \alpha - d_\mu (\pi(\alpha, \beta)) \ .
$$
\end{lemma}

\begin{remark}{\rm{
A bivector $\pi$ is Poisson if and only if $\gamma_\pi =\{\pi,\mu\}$ is primitive
in the sense of Definition~\ref{prim}.
Assume that $\pi$ is a non-degenerate bivector, with inverse $\omega$.
We consider the decomposition \eqref{hll} of the structure
$\mu\in {\mathcal{F}}^{0,1}$,
$$
\mu = \mu_0 + {\mathrm{ad}}_\omega \mu_1 \ ,
$$
where $\mu_0\in {\mathcal{F}}^{0,1}$ and $\mu_1\in {\mathcal{F}}^{1,0}$
are primitive, ${\mathrm{ad}}'_{\pi} \mu_i = 0,\ i=0,1$,
and of weight $1$ and $-1$, respectively.
Then, using Lemma \ref{weight}, we obtain
$$
\gamma_\pi  = \{\pi,\mu\} = \{\pi, \mu_0 + {\mathrm{ad}}_\omega \mu_1 \} = -{\mathrm{ad}}'_{\pi}\mu_0 - {\mathrm{ad}}'_{\pi}{\mathrm{ad}}_\omega \mu_1 = \mu_1 \ .
$$
Thus $\mu=\mu_0 + {\rm{ad}}_\omega \gamma_\pi $, where $\mu_0$ and $\gamma_\pi $ are primitive.
}}
\end{remark}

\subsection{Dualization and composition}
We now dualize the construction of Section \ref{section4.1}.
Let $(A^{*} ,\gamma)$ be a Lie algebroid.
If $\omega \in {\rm{\Gamma}} (\wedge^{2} A^{*})$, then ${\widetilde\mu
  = \{ \gamma ,\omega\}}$ is of shifted bidegree $(0,1)$ and
$\widetilde\mu$ is a Lie algebroid structure on $A$ if and only if
\begin{equation}
\{ [\omega,\omega]_{\gamma} ,\gamma\} =0 \ ,
\end{equation}
while a sufficient condition is
\begin{equation}\label{omega}
{[\omega ,\omega]_{\gamma}=0}\ .
\end{equation}

We shall now combine the two preceding constructions and consider the
following scheme,
$$
\boxed{(A,\mu) \overset{(\pi)}\rightsquigarrow (A^{*},\gamma_\pi   = \{ \pi,\mu\})
\overset{(\omega)} \rightsquigarrow 
 (A , \widetilde\mu   = \{ \gamma_\pi  ,\omega \} = \{\{ \pi, \mu\}, \omega
 \}) \ .}
$$

The following definition is due to Vaisman \cite{V1}.

\begin{definition}
A 2-form satisfying \eqref{omega} when $\gamma = \gamma_\pi = \{\pi, \mu \}$
is called a {\emph {complementary
  $2$-form}} for $\pi$.
\end{definition}

Since, in this case,
$\gamma = \{ \pi,\mu\}$, by Lemma \ref{lemma1}, $[\omega
,\omega]_{\gamma} = \{\{\omega, \gamma\}, \omega\}$ is equal to
$[\omega ,\omega]_{\pi}$, where $[~,~]_\pi$ is the Koszul bracket .

\medskip

Let $\pi$ be an arbitrary bivector and $\omega$ an arbitrary 2-form.
Let us determine sufficient conditions for $\widetilde\mu
= \{\gamma_\pi, \omega\}$
to be a Lie algebroid structure on $A$, {\it i.e.}, to satisfy
\begin{equation}\label{tildemu}
\{ \widetilde\mu ,\widetilde\mu \} =0 \ .
\end{equation}

\begin{proposition}\label{prop3}
(i) Let $\pi$ be a bivector on $(A,\mu)$ such that $\gamma_\pi  = \{\pi,\mu\}$
satisfies $\{\gamma_\pi , \gamma_\pi \} =0$.
A necessary and sufficient condition for $\widetilde\mu = \{ \gamma_\pi , \omega \} =
\{ \{\pi, \mu\}, \omega\}$ to be a Lie
algebroid structure on $A$ is $\{[\omega , \omega ]_\pi , \gamma_\pi  \} =0$.

(ii) Let $\pi$ be a bivector on $(A,\mu)$.
A sufficient condition for $\widetilde\mu = \{ \gamma_\pi , \omega \} =
\{ \{\pi, \mu\}, \omega\}$ to be a Lie
algebroid structure on $A$ is
$$
\left\lbrace
\begin{array}{ll}
[\pi,\pi]_{\mu} = 0 & (\pi \text{ is Poisson}) \ , \\

[\omega ,\omega]_{\pi} = 0 & (\omega \text{ is a complementary
  2-form for} \, \, \pi) \ .
\end{array}
\right.
$$
\end{proposition}

\noindent{\it {Proof}}
Using the Jacobi identity we compute
$$
\{ \widetilde\mu , \widetilde\mu \} = 
\{ \{ \gamma_\pi  , \omega\} , \{ \gamma_\pi ,\omega \}\} =
\{\gamma_\pi , \{ \omega , \{ \gamma_\pi , \omega\} \} \} +
\{ \{ \gamma_\pi , \{ \gamma_\pi  , \omega \}\}, \omega\}
$$
$$
= \{\gamma_\pi , \{ \omega , \{ \gamma_\pi , \omega\} \} \} +
\frac{1}{2}
\{\{\{ \gamma_\pi  , \gamma_\pi \}, \omega \}, \omega \} .
$$

Let us assume that $\gamma_\pi  = \{ \pi,\mu\}$, and
that $\pi$ satisfies $\{\gamma_\pi  ,\gamma_\pi  \} = 0$, which is 
equivalent to \eqref{1}.
Condition \eqref{tildemu} becomes
$$
\{\gamma_\pi , \{ \{\omega ,\gamma_\pi  \} , \omega \} \} =0 \ ,
$$
{\it i.e.},
\begin{equation}\label{5}
\{ \gamma_\pi ,  [\omega ,\omega ]_{\pi} \}=0   \ .
\end{equation}
This proves part (i), and part (ii) follows immediately.
\hfill $\square$

\subsection{Lie algebroid structure defined by a complementary
  $2$-form}
Let us determine an explicit expression for the anchor and
bracket of $(A,\widetilde\mu)$. By the Jacobi identity,
\begin{equation}\label{mutilde}
{ \widetilde\mu = \{ \{\pi,\mu\}, \omega\} = \mu_{1} +\mu_{2}} \ ,
\end{equation}
where we have set
$$
 \mu_{1} = \{\{ \pi, \omega \} ,\mu
\}
\quad\quad {\mathrm {and}}\quad\quad
\mu_{2}  = \{ \pi, \{ \mu , \omega  \}\}
\ .
$$

We set
$$
N = \{\pi, \omega\} \ ,
$$
then $\mu_{1} = \{ N ,\mu\}$ and
we write $\mu_{2} = \{ \pi , \{ \mu, \omega \}\} = \{ \psi, \pi\}$, where $\psi$
is the $3$-form $- \{\mu, \omega \} = - d_{\mu} \omega$.

By definition, the anchor of $(A,\widetilde \mu)$
is $\widetilde \rho$ such that
$\widetilde \rho(X) f =\{ \{ X, \widetilde \mu\},f\} $,
for all $X \in {\rm{\Gamma}} A$ and $f \in C^\infty(M)$.
Then $\widetilde \rho = \rho \circ N$,
where $\rho$ is the anchor of $A$.
In fact, for $X\in
{\rm{\Gamma}} A$ and $f\in C^{\infty}(M)$,
$$
\{\{ X , \mu_{1} \} , f\} = \{\{ X , \{ N ,\mu \}\} ,f\} = \rho
(NX) \cdot f 
$$
and
$$
\{\{ X , \mu_{2}\} , f\} = \{\{ X , \{ \psi, \pi \}\} , f\} = 0\ .
$$

\medskip

Let us consider the bracket defined,
for $X$, $Y \in {\rm{\Gamma}} A$, by
$$
[X,Y]_{\widetilde\mu} = \{\{ X , \widetilde\mu \},Y\}\ .
$$
By Lemma \ref{bracketN}, the bracket
$[~,~]_{\mu_{1}}$ is the
bracket $[~,~]^\mu_{N}$ recalled in \eqref{deformed}.
The theorem below follows from Lemmas \ref{lemma2},
\ref{bracketN} and \ref{lemma3}.
\begin{theorem}\label{theorem3} 
Let $\pi$ be a bivector and $\omega$ a 2-form on
$(A,\mu)$. Then $\pi$ and $\omega$ satisfy \eqref{tildemu}
if and only if the
bracket of sections of $A$ defined by
\begin{equation}\label{6}
[X,Y]_{\widetilde\mu} = [X,Y]^\mu_{N} - \pi^{\sharp} (i_{X\wedge Y} d_{\mu}
\omega) \ ,
\end{equation}
for all $X$, $Y\in {\rm{\Gamma}} A$, where $N = \pi^{\sharp}\circ
\omega^{\flat}$,
is a Lie algebroid bracket with anchor $\rho \circ N$.
\end{theorem}

In order to compare \eqref{6} with Formula (3.3) in \cite{V1}, we
remark that $B = - N$, so that $[~,~]'_E$ is the opposite of
$[~,~]_{\widetilde \mu}$. As a corollary of Proposition \ref{prop3} (ii)
and Theorem \ref{theorem3}, we obtain the following results which were proved
in \cite{V1}. 

\begin{corollary} If $\pi$ is a Poisson bivector and $\omega$
is a complementary 2-form for $\pi$, then

(i) Formula \eqref{6} defines a Lie bracket on the space of sections of $A$, and

(ii) if, in addition, $d_\mu \omega = 0$, then bracket $[~,~]^\mu_{N}$ , where
$N = \pi^{\sharp} \circ \omega^{\flat}$, is a Lie bracket.
\end{corollary}

In part (ii) of this corollary the assumption that $\omega$ be
$d_\mu$-closed can be replaced by the weaker assumption that, 
for all $X$ and $Y \in {\rm{\Gamma}} A$, $i_{X\wedge Y} d_\mu
\omega \in {\rm{ker}}(\pi^\sharp)$.

\begin{remark}{\rm{If $\pi$ is a non-degenerate Poisson bivector, its
inverse $\omega$ is
a complementary 2-form for $\pi$. In fact \cite{KM},
$[\omega , \cdot \,  ]_{\pi}=d_{\mu}$ and therefore $[\omega,\omega]_{\pi} =
d_{\mu}\omega = 0$. In this case $N = {\rm{Id}}_A$ and $\mu_{2} = 0$,
therefore $\widetilde\mu =\mu$.}}
\end{remark}

\subsection{The modular class of $(A,\widetilde \mu)$}\label{modclass}

Consider a Poisson bivector $\pi$ on $A$ and
a complementary $2$-form $\omega$ with respect to $\pi$.
Assume that $A$ is orientable and let
$\lambda$ be a nowhere-vanishing section of $\wedge^{\rm{top}}(A^*)$
that defines an isomorphism, $*_\lambda$, from multivectors to forms.
Let $d^\pi_\omega = - [\omega, . ]_\pi$ be the Lie algebroid cohomology
operator of $A$ with structure $\widetilde \mu = \{\{\pi,
\mu\},\omega\}$.
Each of the operators on the sections of $\wedge^\bullet A$,
$$
\partial^\pi_\omega = [d_\pi , i_\omega] \ ,
$$
and
$$
\partial^\pi_{\omega,\lambda} = - (*_\lambda)^{-1} d^\pi_\omega \, *_\lambda
$$
generates $[~,~]_{\widetilde \mu}$
and also has square $0$ since $\omega$ is a complementary $2$-form with respect to $\pi$.
The $1$-form  $\xi_{\pi,\omega,\lambda}$ on $A$ defined by
$$
\partial^\pi_{\omega,\lambda} - \partial^\pi_\omega =
i_{\xi_{\pi,\omega,\lambda}}
$$
is a $d_{\widetilde \mu}$-cocycle. Its class is the modular class of
the Lie algebroid $(A,\widetilde \mu)$ \cite{V4}.

\medskip

In the following sections, $A$ denotes a vector bundle over a
manifold $M$, and we let $(A,\mu)$ be a Lie algebroid,
so that,
by assumption, $\{\mu,\mu\}=0$. We will sometimes abbreviate $(A,\mu)$ by $A$.

\section{What is a $PN$-structure on a Lie algebroid?}\label{sectionPN}

We have reviewed the Nijenhuis structures in Section \ref{sectionNij}. We now
consider Nijenhuis structures on Lie algebroids equipped with a Poisson
structure.

\subsection{Compatibility}
Given a bivector $\pi$ and a $(1,1)$-tensor $N$
on $(A,\mu)$, we consider both
$$
\mu_N
= \{N, \mu\} \ ,
$$
which defines an anchor $\rho \circ N$
and a bracket  $[~,~]^\mu_N$ on $A$,
and
$$
\gamma_\pi = \{\pi,\mu\} \ ,
$$
which defines an anchor
$\rho \circ \pi^\sharp$ and a bracket  on $A^*$
that we have denoted by $[~,~]_\pi$. We assume that
\begin{equation}\label{skew}
N \circ \pi^\sharp = \pi^\sharp \circ N^* \ ,
\end{equation}
where $N^*$ denotes the transpose of $N$ satisfying $<N(X),\alpha> =
<X,N^*(\alpha)>$, for all $X \in {\mathrm{\Gamma}} A$ and $\alpha \in {\mathrm{\Gamma}}(A^*)$, 
so that $N \circ \pi^\sharp$ defines a bivector $\pi_N$ by
$\pi_N^\sharp = N \circ
\pi^\sharp$. Then,
$$
\pi_N = \frac{1}{2} \{\pi, N\} \ .
$$
We introduce a \emph{compatibility condition} for $\pi$ and $N$ by requiring
that the bracket $[~,~]^\mu_N$ twisted by $\pi$, which is $\{\pi,\{N,
\mu\}\}$,
be equal to the bracket
$[~,~]_\pi$ deformed by $N^*$, which is $\{\{\pi,
\mu\},N\}$.
Thus we set
\begin{equation}\label{comp}
C_\mu(\pi,N) = \{\pi,\{N, \mu\}\} + \{N,\{\pi, \mu\}\} \ ,
\end{equation}
which is a section of
$\wedge^2 A \otimes A^*$.

\begin{definition}
A bivector $\pi$ and a $(1,1)$-tensor $N$
on $(A,\mu)$ are called \emph{compatible} if they satisfy \eqref{skew} and
$$
C_\mu(\pi,N)= 0 \ .
$$
A {\emph{PN-structure}} on $(A,\mu)$
is defined by a Poisson bivector and a Nijenhuis
tensor on $(A,\mu)$ which are compatible.
\end{definition}

\subsection{$PN$ structures and Lie bialgebroid structures}
A necessary and sufficient condition
for $(\mu_N,\gamma_\pi)$ to define a Lie bialgebroid
structure on $(A,A^*)$ is $\{\mu_N + \gamma_\pi, \mu_N + \gamma_\pi\}
= 0$. When $N$ is a Nijenhuis $(1,1)$-tensor, $\{\mu_N,\mu_N\}=0$, and
when $\pi$ is a Poisson bivector,
$\{\gamma_\pi,\gamma_\pi\}=0$. Therefore in this case the condition
$\{\mu_N + \gamma_\pi, \mu_N + \gamma_\pi\} = 0$ is
equivalent to $\{\mu_N,\gamma_\pi\}=0$.

\begin{lemma} Let $C'_\mu(\pi,N) = 2 \{\mu_N,\gamma_\pi\}$. Then
  $C'_\mu(\pi,N) =\{\mu,C_\mu(\pi,N)\}$.
\end{lemma}

\noindent{\it{Proof}}
By the Jacobi identity,
$$
C'_\mu(\pi,N) = 2 \{\{ N,\mu\}, \{\pi,\mu\}\} =
\{\{ N,\mu\}, \{\pi,\mu\}\}
+ \{\{ \pi,\mu\}, \{N,\mu\}\}
$$
$$
 \quad \quad \quad \quad \quad \quad=  \{\{\{N,\mu\}, \pi\},\mu\} +
\{\{\{\pi,\mu\}, N\}, \mu\} = \{\mu, C_\mu(\pi,N)\} \ . \quad \quad
\quad   \quad \quad   \quad\square
$$

\begin{theorem} Let $N$ be a Nijenhuis $(1,1)$-tensor and $\pi$ a
Poisson bivector on $(A,\mu)$.

(i) The vanishing of $\{\mu, C_\mu(\pi,N)\}$ is a
necessary and sufficient condition for $(\mu_N,\gamma_\pi)$
to define a Lie bialgebroid structure on $(A,A^*)$.
In particular, if $\pi$ and $N$ are compatible, then
$(\mu_N,\gamma_\pi)$ is a Lie bialgebroid structure.

(ii) If the $d_\mu$-exact $1$-forms generate ${\rm{\Gamma}}(A^*)$ locally as a
$C^\infty(M)$-module, then a Poisson bivector $\pi$ and a Nijenhuis
tensor $N$ define a $PN$-structure on $(A,\mu)$ if and only if
the pair $(\mu_N,\gamma_\pi)$
defines a Lie bialgebroid structure on $(A,A^*)$.
\end{theorem}

\noindent{\it{Proof}}
Only (ii) needs to be proved.
From $C'_\mu(\pi,N) =\{\mu,C_\mu(\pi,N)\}$ we obtain
$$
\{C_\mu(\pi,N), \{\mu,f\}\} = \{\{C_\mu(\pi,N),\mu\},f\} =
- \{C'_\mu(\pi,N),f\} \ ,
$$
for all $f \in C^\infty(M)$. Thus $C'_\mu(\pi,N) = 0$ implies
$C_\mu(\pi,N)(d_\mu f , \cdot )= 0$, for all $f \in C^\infty(M)$.
Under the assumptions of part (ii) of the theorem, 
$C_\mu(\pi,N)$ vanishes identically
since it is $C^\infty(M)$-linear.
\hfill $\square$

\medskip

The equivalence stated
in the theorem was proved in \cite{yks96} for the case when $A$ is 
the tangent bundle of a manifold.
This equivalence may fail for Lie algebroids
which are not tangent bundles, a fact observed by Grabowski and
Urbanski \cite{GU}. 

\begin{remark} {\rm{The compatibility condition
$C_\mu(\pi,N)= 0$ implies that the brackets
$[~,~]_{\pi_N}$, $([~,~]^\mu_N)_\pi$ and $([~,~]_\pi)_{N^*}$
coincide. In fact, from $\{\pi,\{N,\mu\}\} = \{\{\pi,\mu\}, N\}$, we
obtain
$\{\mu,\{N,\pi\}\} = \{\{\mu,N\}, \pi\} + \{N,\{\mu,\pi\}\} =
2 \{\{\mu,N\}, \pi\} = 2 \{N, \{\mu, \pi\}\}$ 
or $\{\mu , \pi_N\} = \{\mu_N, \pi\} = \{\{\pi, \mu\},N \}$.}}
\end{remark}

\begin{remark}\label{remark7}
{\rm{The condition $C'_\mu(\pi,N)= 0$ is
      equivalent to each of the following:

\noindent $\bullet$ The operator $d_N = \{\{N,\mu\}, \cdot \} = [i_N,d_\mu]$ is a
derivation of $[~,~ ]_\pi$,

\noindent $\bullet$
$({\rm{\Gamma}}(\wedge^\bullet A^*), [~,~]_\pi, d_N)$ is a
differential Gerstenhaber algebra,

\noindent $\bullet$ The operator $d_\pi = \{\{\pi,\mu\}, \cdot \} =
[\pi,\cdot]_\mu$ is a derivation of $[~,~]^\mu_N$,

\noindent $\bullet$
$({\rm{\Gamma}}(\wedge^\bullet A), [~,~]_N^\mu, d_\pi
)$ is a
differential Gerstenhaber algebra.}}
\end{remark}

\section{On compatible Poisson structures}\label{sectioncompatible}
For some of the results in this section, see \cite{D} and earlier articles cited there.

Two Poisson bivectors on $(A,\mu)$ are said to be {\emph{compatible}},
or to form a {\emph{Hamiltonian pair}} or to define 
a {\emph {bi-Hamiltonian structure}}, if their sum is a Poisson bivector.
Thus Poisson bivectors $\pi$ and $\pi_1$ are compatible if and only if
$\{\{\pi,\mu\},\pi_1\} = 0$.

Let $\pi$ be a bivector on $(A,\mu)$ and $N$ a $(1,1)$-tensor.
Assume that $N \circ \pi^\sharp$ is a bivector, {\it i.e.},
$N \circ \pi^\sharp
=\pi^\sharp  \circ N^*$.  We have set
 $\pi_N^\sharp = N \circ \pi^\sharp$. Then
$$
\pi_N = \frac{1}{2} \{\pi,N\} \ .
$$
In particular, if a bivector $\pi$ is non-degenerate
and has inverse $\omega$, then $N = \{\pi_N,\omega\}$.

\begin{proposition}\label{compatible}
Assume that $\pi$ is a non-degenerate Poisson bivector on $(A,\mu)$
with inverse $\omega$, and
$N$ is a $(1,1)$-tensor such that $\pi$ and $\pi_N$
satisfy $\{\{\pi,\mu\},\pi_N\} = 0$. Then
\begin{equation}\label{C}
\{\{\mu,\pi\},N\} + \{\{\mu,N \},\pi \} =0 \ ,
\end{equation}
and
\begin{equation}\label{C0}
\{\{N,\mu\},\omega\} = 0 \ .
\end{equation}
If, in addition, $\pi_N$ is a Poisson bivector, then
\begin{equation}\label{C1}
\{\{\mu,\pi_N\},N\} + \{\{\mu,N \},\pi_N \} =0 \ .
\end{equation}
\end{proposition}

\noindent{\it{Proof}}
Since $\omega$ is the inverse of $\pi$, $\{\pi, \omega \}$ is the
identity of $A$ and, by \eqref{inverse}, for $u \in {\mathcal F}^{p,q}$,
\begin{equation}\label{id}
\{u, \{\pi,\omega \}\} = (p - q) u \ .
\end{equation}
By assumption, $\{\{\pi,\mu\},\pi\} = 0$
and $\{\{\pi,\mu\},\pi_N\} = 0$.

(i) From the compatibility of $\pi$ and $\pi_N$ and the Jacobi identity,
we derive
$$
0= \{\{\{\pi_N,\mu\},\pi\},\omega\}
=\{\{\{\pi_N,\mu\},\omega\},\pi\} + \{\{\pi_N,\mu\}, \{\pi,\omega\}\}
\ .
$$
Because $\pi$ is a Poisson bivector,
$\{\mu, \omega\}=0$. Therefore
$$
0
=\{\{\{\pi_N,\omega\},\mu\},\pi\} + \{\{\pi_N,\mu\}, \{\pi,\omega\}\}
\ .
$$
Using the relations
$\pi_N = \frac{1}{2} \{\pi,N\}$, $N=\{\pi_N, \omega \}$, and
\eqref{id}, we obtain
$$
0 = \{\{N,\mu\},\pi\} + \frac{1}{2} \{\{\pi,N\},\mu\}
= \{\{N,\mu\},\pi\} + \frac{1}{2} \{\pi , \{N,\mu\}\} +
\frac{1}{2} \{\{\pi,\mu\},N\} 
$$
$$
= \frac{1}{2} \{\{N, \mu\},\pi\} + \frac{1}{2} \{\{\pi,\mu\},N\} \ .
$$

(ii) To prove \eqref{C0}, we compute
$$
\{\{\pi, \{N,\mu\}\},\pi\} = \{\{\{\pi,N\},\mu\},\pi\} +
 \{\{N, \{\pi,\mu\}\},\pi\} =
4 \{\pi_N,\{\mu,\pi\}\} = 0 \ ,
$$
and we know that the
vanishing of $\{\{\pi, \{N,\mu\}\},\pi\}$ is equivalent to
$\{\{N,\mu\}, \omega \} =0$.

(iii) To prove \eqref{C1},
we use the assumption $\{\{\pi_N, \mu\},\pi_N\} = 0$ to
obtain
\begin{align*}
0 & =  \{\{\{\pi_N, \mu\},\pi_N\}, \omega\} =
\{\{\pi_N,\mu\},N\} +
\{\{\{\pi_N,\mu\},\omega\}, \pi_N\} \\
& = \{\{\pi_N,\mu\},N\} +
\{\{N,\mu\}, \pi_N\} \ ,
\end{align*}
thus proving \eqref{C1}. \hfill $\square$

\medskip

From \eqref{C}, it follows that
\begin{equation}\label{eq4}
\{\pi_N,\mu\} = \frac{1}{2} \{\{\pi,N\},\mu\} = \{\{\pi,\mu\},N\} =  \{\{\mu,N\},\pi\} \ ,
\end{equation}
and from \eqref{C1}, it follows that
\begin{equation}\label{eq5}
\frac{1}{2} \{\{\pi_N,N\},\mu\} = - \{\{\mu,\pi_N\},N\} =
\{\{\mu,N\},\pi_N\} \ . 
\end{equation} 

\begin{lemma}\label{Nsquare}
When $N = \{\sigma, \tau\}$, where $\sigma$
is a bivector and $\tau$ a $2$-form, then
$$
N^2 = - \frac{1}{2} \{\{N,\sigma\},\tau\} \ .
$$
\end{lemma}

\noindent{\it Proof} This formula is proved by a simple
calculation. \hfill $\square$

\medskip

The following essential result in the theory of bi-hamiltonian systems
was proved by Magri and Morosi in \cite{MM} and
also by Gelfand and Dorfman \cite{D} in the algebric framework of
Hamiltonian pairs and by Fuchssteiner and Fokas \cite{FF} 
in their study of Hamiltonian structures for evolution equations.
See \cite{KM} for the case of Lie algebroids.

\begin{theorem}\label{Nij}
Let $\pi$ and
$\pi_1$ be compatible Poisson
structures, with $\pi$ non-degenerate. Set $N= \pi_1^\sharp \circ
(\pi^\sharp)^{-1}$. Then

(i) the Nijenhuis torsion of
the $(1,1)$-tensor $N$ vanishes.

(ii) the pair $(\pi,N)$ is a PN-structure,

(iii) the pair $(\pi_1,N)$ is a PN-structure.
\end{theorem}

\noindent{\it{Proof}}
(i) Let $\omega$ be the inverse of $\pi$. Then $N = \{\pi_1,\omega\}$.
Applying Lemma \ref{Nsquare} to $\sigma=\pi_1$ and $\tau = \omega$ ,
we obtain
from Formula \eqref{a}:
\begin{align*}
    2 \, {\mathcal T}_\mu N &= \{N, \{N,\mu\} \}- \{N^2, \mu\}\\
    &=  \{N, \{N,\mu\} \} +
    \frac{1}{2}\{\{\{N,\pi_1\},\omega\}, \mu\} \ .
\end{align*}
Using \eqref{C0} and \eqref{C1} yields
\begin{align*}
    2 \, {\mathcal T}_\mu N &= \{N, \{N,\mu\} \} + \{\{\{N,\mu\}.\pi_1\}, \omega\}\\
    &=  \{N,\{N,\mu\}\} + \{\{N,\mu\},N\}  = 0 \ .
\end{align*}

(ii) Equation \eqref{C} expresses the vanishing of $C_\mu(\pi,N)$.

(iii) Equation \eqref{C1} expresses the vanishing of $C_\mu(\pi_1,N)$.
\hfill $\square$

\medskip

In \cite{KM} and in the references cited above, 
it is proved more generally that
each $\pi_k$ defined by $\pi_k^\sharp = N^k \circ
\pi^\sharp$, $k \in \mathbb N$,
is a Poisson bivector, and the bivectors $\pi_k$ are pairwise compatible.
The preceding theorem admits a converse.

\begin{theorem}\label{converse}
If $\pi$ and $\pi_1$ are non-degenerate Poisson bivectors, and
if  $N= \pi_1^\sharp \circ
(\pi^\sharp)^{-1}$ has vanishing Nijenhuis torsion, then $\pi$ and
$\pi_1$ are compatible. 
\end{theorem}

\noindent{\it{Proof}}
Assume that $\pi$ is non-degenerate and let $\omega$ be its inverse.
Since $\pi$ is a Poisson bivector, $\{\mu, \omega\}=0$.
From this fact, the fact that $\pi_1$ is Poisson and the formula
$N = \{\pi_1,\omega\}$, we obtain 
\eqref{C1} which implies
\begin{equation}\label{C1'}
\{\{N,\pi_1\},\mu\} = 2 \{\{N,\mu\},\pi_1\} \ .
\end{equation}
From Lemma \ref{Nsquare} applied to $N = \{\pi_1,\omega\}$, we obtain
$$
2 {\mathcal T}_\mu N = \{N,\{N,\mu\}\} + \frac{1}{2} \{\{\{N,
\pi_1\},\omega\},\mu\} \ .
$$
Using $\{\mu,\omega\} = 0$  and \eqref{C1'}, we obtain
$$
2 {\mathcal T}_\mu N = \{N,\{N,\mu\}\} + \{\{\{N, \mu\},\pi_1\},\omega\} \ ,
$$
whence
\begin{equation}\label{C1''}
{\mathcal T}_\mu N = \frac{1}{2} \{\{\{N, \mu\},\omega\},\pi_1\} \ .
\end{equation}
We now assume that $\pi_1$ also is non-degenerate, and we denote its
inverse by $\omega_1$. Then $2\{{\mathcal T}_\mu N, \omega_1\} = \{\{N,
\mu\},\omega\}$ and the vanishing of ${\mathcal T}_\mu N$ implies the
vanishing of $\{\{N, \mu\},\omega\}$.
We then remark that the vanishing of  $\{\{N, \mu\},\omega\}$ is
equivalent to the vanishing of $\{\{\pi, \{N,\mu\}\}, \pi\}$.
Since $\{\{\pi_1, \mu\}, \pi\} = \frac{1}{4} \{\{\pi,
\{N,\mu\}\}, \pi\}$, the vanishing of  ${\mathcal T}_\mu N$ implies
that $\pi$ and $\pi_1$ are compatible.
\hfill $\square$

\medskip

\begin{remark}{\rm{Let $\pi$ be a non-degenerate Poisson bivector with
inverse $\omega$ and let $N$ be a $(1,1)$-tensor. Assume that 
$\pi_N$ defined by $(\pi_N)^\sharp = N \circ \pi^\sharp$ is a Poisson
bivector. In view of Lemma \ref{lemma3} and Lemma \ref{lemma11} below,
Formula \eqref{C1''} means that 
$$
({\mathcal T}_\mu N)(X,Y) = \frac{1}{2} (\pi_N)^\sharp(i_{X \wedge Y}
  (d_N\omega)) \ ,
$$
for all $X$ and $Y \in {\rm{\Gamma}}A$.}}
\end{remark}

In the next sections we shall review and compare
the $P\Omega$- and $\Omega N$-structures of Magri
and Morosi \cite{MM} and the Hitchin pairs of Crainic \cite{C}.

\section{What is a $P\Omega$-structure on a Lie algebroid?}\label{sectionpomega}

In \cite{MM}, Magri and Morosi defined the $P \Omega$- and
$\Omega N$-structures on manifolds and, more recently, in his study of
generalized complex structures, Crainic
defined Hitchin pairs on manifolds \cite{C}.
These notions admit straightforward generalizations to 
the case of Lie algebroids which we now define and study. 
When the Lie algebroid is $TM$ with its standard Lie algebroid
structure, these definitions recover the classical case.
Most of the results in this section are particular cases of the general
theorem of Antunes\footnote{The convention for the definition of
  $\omega^\flat$ in \cite{A} is the opposite of ours.} 
on Poisson quasi-Nijenhuis structures with
background, see \cite{A}, Theorem 4.1.

If $N$ is a $(1,1)$-tensor on a Lie algebroid $(A,\mu)$, 
let $\mu_N= \{N,\mu\}$ be the deformed bracket satisfying
\eqref{deformed}. Define the operator on forms $i_N = \{N, \cdot\}$
and let $d_N$ be the operator considered in Remark \ref{remark7},
$$
d_N = [i_N, d_\mu] \ ,
$$
where $[~,~]$ is the graded commutator.
In particular, if a form $\alpha$ is
 $d_\mu$-closed, then
 $d_\mu(\alpha_N) = - d_{N}\alpha$, where $\alpha_N = i_N \alpha$.
The following simple lemma was proved in \cite{KM}.
We present an alternative proof.

\begin{lemma}\label{lemma11}
Let $N$ be a $(1,1)$-tensor on $(A, \mu)$.
The operators on forms $d_N= [i_N, d_\mu]$ and $d_{\mu_N} = \{ \mu_N,
 \cdot \}$ coincide. 
\end{lemma}

\noindent{\it{Proof}}
For any form $\alpha$,
$d_{\mu_N} \alpha = \{ \mu_N , \alpha \} = \{ \{ N, \mu \}, \alpha \}
= \{ N, \{ \mu, \alpha \} \} - \{\mu, \{ N, \alpha \} \}$, while
$d_N \alpha = [i_N, d_\mu](\alpha) =i_N \{\mu, \alpha \} - \{ \mu, i_N
\alpha \} =  \{ N, \{ \mu, \alpha \} \} - \{\mu, \{ N, \alpha \} \}$.
\hfill $\square$

\medskip

Let $N$ be a $(1,1)$-tensor and $\omega$ a $2$-form on $(A, \mu)$
such that 
\begin{equation}\label{c}
 \omega ^\flat  \circ N  = N^* \circ \omega^\flat \ .
\end{equation}
Then $\omega_N$ defined by $\omega_N^\flat = \omega^\flat \circ N$ is 
a $2$-form and
$$
\omega_N = \frac{1}{2} i_N \omega = \frac{1}{2} \{N, \omega\} \ .
$$ 

Let $\pi$ be a bivector and let $\omega$ be a $2$-form on the Lie
algebroid $(A,\mu)$.
Set $N =  \pi^\sharp \circ \omega^\flat$. Then 
\eqref{c} is satisfied and
$$
N = \{\pi,\omega\} \ .
$$

We shall now prove identities relating $\pi$, $\omega$ and $N$ when $N =
\{\pi, \omega\}$.
\begin{lemma}\label{lemmad_N}
Let $\pi$ be a bivector and $\omega$ a $2$-form on $(A, \mu)$,
and let $N$ be the $(1,1)$-tensor $N = \pi^\sharp \circ \omega^\flat$.
The $2$-form $\omega$ satisfies
\begin{equation}\label{d_N}
\frac{1}{2} [\omega, \omega]_\pi + d_N \omega = - \frac{1}{2}d_\mu(\omega_N)  \ .
\end{equation}
\end{lemma}
\noindent{\it{Proof}} In fact, by the Jacobi identity,
\begin{align*}
 [\omega, \omega]_\pi & = 
\left \{\{ \omega,  \{ \pi, \mu  \} \} , \omega  \right \} = 
\{\{\{ \omega,  \pi \}, \mu  \} , \omega  \} + 
\{\{\pi, \{\omega, \mu \} \} , \omega  \} \\
& = \{\{\omega,\pi\},\{\mu, \omega \}\}
+ \{\{\{\omega,\pi\}, \omega\},\mu\}
+\{\{\pi,\omega\},\{\omega, \mu \}\} \\
&
= 2 \{\{\omega,\pi\},\{\mu, \omega \}\} -
\{ \{\{ \pi, \omega \}, \omega\}, \mu\} \ .
\end{align*}
Since $N = \{\pi, \omega\}$ and $d_\mu \omega = \{\mu, \omega\}$, 
we obtain $ [\omega, \omega]_\pi = - 2 \{N,\{\mu,\omega\}\} +
\{\mu,\{N,\omega\}\}$, hence \eqref{d_N}.
\hfill $\square$

\medskip

\begin{remark}{\rm{
When $(\pi,N)$ is
a $PN$-structure, $({\rm{\Gamma}}(\wedge^\bullet(A^*)), [~,~]_\pi,
d_N)$ is a differential graded Lie algebra, and this formula expresses
the fact that $2$-forms  such that $i_N\omega$ is $d_\mu$-closed 
are Maurer-Cartan elements in this DGLA.}}
\end{remark}

We shall prove two additional identities relating $\pi$, $\omega$ and $N$.
We recall that $[\pi,\pi]_\mu = \{\{\pi,\mu\},\pi\}$ and that
$C_\mu(\pi,N)$ is defined by Formula \eqref{comp}.
\begin{lemma}\label{lemmarelations}
A bivector $\pi$, a $2$-form $\omega$ and a $(1,1)$-tensor 
$N$ related by  
$N = \pi^\sharp \circ \omega^\flat$ satisfy the relations
\begin{equation}\label{C=0}
\{[\pi,\pi]_\mu, \omega\} = \{\{\pi, d_\mu \omega\}, \pi\} + C_\mu(\pi,N)
\ ,
\end{equation}
and
\begin{equation}\label{T=0}
\{\{[\pi,\pi]_\mu, \omega\}, \omega\} = \{\{\pi, d_\mu \omega\},
\pi\},\omega\} - \{\{\pi,N\},d_\mu \omega\} + 2 \{\pi,d_N\omega\} 
+ 4 {\mathcal T}_\mu N \ .
\end{equation}

\end{lemma}
\noindent{\it{Proof}} Applying the
Jacobi identity, we obtain
$$
\begin{array}{lll}
\{\{\{\pi,\mu\},\pi\},\omega\} & = &
\{\{\{\pi,\mu\},\omega\},\pi\}  
+ \{\{\pi, \mu\},\{\pi,\omega\}\}\\
& = &\{\{\pi, \{\mu , \omega\}\},\pi\} +  \{\{\{\pi,\omega\},\mu\},\pi\}  
+ \{\{\pi, \mu\},\{\pi,\omega\}\} \ .
\end{array}
$$
Therefore, for $N= \{\pi,\omega\}$, we obtain \eqref{C=0}.
Furthermore,
\begin{align*}
& \{C_\mu(\pi,N),\omega\} 
= \{\omega, \{\{N, \mu\},\pi\}\} +
\{\omega, \{\{\pi, \mu\},N\}\} \\
& = \{\{N,\mu\},\{\omega, \pi\}\} + \{\pi, \{\{N,\mu\},\omega\}\} 
+ \{ \omega, \{\pi , \{ \mu,N\} \} \}  + \{\omega,  \{\{\pi, N\},\mu\}\} \\
& = - \{\{N,\mu\},N\} + \{\pi, d_N\omega\} + 
\{ \{ \omega, \pi \},\{ \mu,N\} \}\}  \\
& \,  \quad   + \{ \{ \{ \mu,N \}, \omega\}, \pi\}  
+ \{ \{\omega,\{\pi, N\}\},\mu\} + \{\{\pi,N\},\{\omega, \mu\}\} \\
&
= 2 \{N,\{N,\mu\}\} + 2 \{\pi, d_N\omega\} - \{\mu, \{\{N,\pi\},\omega \}\}  
+
\{\{N,\pi\},d_\mu \omega \}\ .
\end{align*}
We shall now make use of Formula \eqref{a} and 
Lemma \ref{Nsquare} which imply
$$ 4 {\mathcal T}_\mu N = 2 \{N, \{N,\mu\}\} - 
\{\{\{N,\pi\},\omega\}, \mu\} \ .
$$ 
Thus
$$ \{C_\mu(\pi,N),\omega\} 
=  \{\{N,\pi\},d_\mu \omega\} + 2 \{\pi,d_N
\omega\} + 4 {\mathcal T}_\mu N \ .
$$
Whence the relation \eqref{T=0}. \hfill $\square$
\begin{definition}
A bivector $\pi$ and a $2$-form $\omega$ on a Lie algebroid $(A,\mu)$
define a {\it{$P \Omega$-structure}}
if $\pi$ is a Poisson bivector, $\omega$
is $d_\mu$-closed, and
$d_\mu(\omega_N) =0$, where $N = \pi^\sharp \circ \omega^\flat$ and $\omega_N$ is the
$2$-form such that $\omega_N^\flat = \omega^\flat \circ N$.
\end{definition}

Since $d_\mu\omega = 0$, the condition $d_\mu(\omega_N) = 0$ in the
definition of a $P \Omega $-structure is
equivalent to the condition $d_N \omega = 0$, so that
$P\Omega$-structures can be characterized as follows. 
\begin{proposition}
A Poisson bivector $\pi$ and a $d_\mu$-closed $2$-form $\omega$ on $(A,\mu)$
define a $P \Omega$-structure if and only if $d_N \omega = 0$, where
$N = \pi^\sharp \circ \omega^\flat$.
\end{proposition}

Since $N = \{\pi, \omega \}$ 
and $\omega_N = \frac{1}{2} i_N \omega = \frac{1}{2} \{ N, \omega\}$,
in terms of the big bracket,
the conditions in the definition  of $P \Omega$ structure are
$$\{\{\pi, \mu\},\pi\} = 0, \quad \quad \{\mu, \omega \}= 0, \quad
\quad {\mathrm{and}} \quad \quad \{\mu,\{N,\omega\}\} =0 \ .
$$

We can relate the notion of $P \Omega$ structure to that of
complementary $2$-form.

\begin{theorem}\label{POmegacompl}
A Poisson bivector $\pi$ and a $2$-form $\omega$ on $(A,\mu)$
define a $P \Omega$-structure if and only if $\omega$ is a
$d_\mu$-closed complementary $2$-form for $\pi$.
\end{theorem}
\noindent{\it{Proof}} In fact, by the Jacobi identity,
$$
[\omega,\omega]_\pi = \{\{\omega, \{\pi, \mu\}\}, \omega\} = 
\{\{\{\omega, \pi\}, \mu\}, \omega\}+ \{\{\{\omega, \mu\},\pi\},
  \omega\}  
$$
$$
= - \{\{N,\mu\}, \omega \} + \{\{\pi, d_\mu \omega\}, \omega \} 
= - d_N\omega + \{N, d_\mu \omega\}  \ .
$$
Thus, if $\omega$ is both $d_\mu$- and $d_N$-closed, it is a
complementary $2$-form. Conversely is $\omega$ is a $d_\mu$-closed
complementary $2$-form, then $d_N\omega =0$.
\hfill $\square$

\medskip

The following theorem generalizes a result of \cite{MM}.

\begin{theorem}\label{theorempomega}
If a Poisson bivector $\pi$ and a $d_\mu$-closed
$2$-form $\omega$ on $(A,\mu)$
define a $P \Omega$-structure, then the pair $(\pi,N)$, where $N =
\pi^\sharp \circ \omega^\flat$, is a $PN$-structure.
Conversely, if $(\pi,N)$ is a $PN$-structure and $\pi$ is
non-degenerate, then $(\pi, \omega)$, where $\omega^\flat =
(\pi^\sharp)^{-1} \circ N$ is a $P\Omega$-structure.
\end{theorem}

\noindent{\it{Proof}}
If  $(\pi, \omega)$ is a $P\Omega$-structure, Equation \eqref{skew} is
obviously satisfied. 
It follows from \eqref{C=0} that, when $\pi$ is a Poisson bivector and
$\omega$ is a
$d_\mu$-closed $2$-form,  
$\pi$ and $N = \{\pi, \omega\}$ are compatible.
It follows from \eqref{T=0} that, 
if in addition $d_N \omega =0$, then 
$ {\mathcal T}_\mu N = 0$. Therefore
$(\pi,N)$ is a $PN$-structure.

Conversely, it is clear that \eqref{skew} implies that $\omega$ is a
$2$-form.
Assume that $\pi$ is non-degenerate 
and let $\tau$ be the $2$-form inverse of $\pi$.
Then $\{ \{\pi, \tau\} , \cdot \} = \bf I$.
Applying $\{\tau, \cdot\}$ 
to both sides of equation \eqref{C=0} yields
$$
\{\{[\pi,\pi]_\mu, \omega\}, \tau\} = \{\{\{\pi, d_\mu \omega\},
\pi\}, \tau\}
 + \{C_\mu(\pi,N),\tau\} \ .
$$
Now, by \eqref{inverse},
$$
 \{\{\{\pi, d_\mu \omega\},\pi\}, \tau\}
= \{\{\pi, d_\mu \omega\},\{\pi , \tau\}\} +
\{\{\{\pi, d_\mu \omega\}, \tau\},\pi\}\}
= 4 \{d_\mu \omega, \pi \} \ .
$$
Applying $\{\tau, \cdot\}$ once more yields $d_\mu \omega =0$.
Therefore, if $\pi$ is a non-degenerate 
Poisson bivector and $\pi$ and $N$ are
compatible, then $\omega$ is $d_\mu$-closed.
Applying $\{\tau, \cdot\}$ to both sides of equation \eqref{T=0} then
yields $d_N \omega = 0$ when $[\pi,\pi]_\mu$, $C_\mu(\pi,N)$ and ${\mathcal
  T}_\mu N$ all vanish.  
\hfill $\square$

\begin{corollary}\label{torsion}
If $\pi$ is a Poisson bivector and $\omega$ is
a closed complementary
$2$-form for $\pi$, then
 the Nijenhuis torsion ${\mathcal T}_{\mu}(N)$
of $N = \pi^{\sharp} \circ
\omega^{\flat}$ vanishes.
\end{corollary}

\noindent{\it {Proof}}
This corollary follows from Theorems \ref{POmegacompl} and \ref{theorempomega}.
\hfill $\square$

\section{What is an $\Omega N$-structure on a Lie algebroid?}\label{8}
\begin{definition}\label{defomegaN}
A $2$-form $\omega$ and a $(1,1)$-tensor $N$ on a Lie algebroid
$(A,\mu)$ define an $\Omega N$-structure if
$\omega^\flat \circ N = N^* \circ \omega^\flat$, $\omega$ is
$d_\mu$-closed, $N$ is a Nijenhuis tensor,
and
$d_\mu(\omega_N) = 0$, where $\omega_N^\flat = \omega^\flat \circ N$.
\end{definition}
Since $\omega_N = \frac{1}{2} \{N, \omega \}$,
the conditions in the definition  of an $\Omega N$ structure, in
addition to $\omega^\flat \circ N = N^* \circ \omega^\flat$, are,
in terms of the big bracket,
$$   \quad \quad \{\mu, \omega \}= 0,  \quad
\quad \{N,\{N,\mu\}\}- \{N^2,\mu\} = 0, \quad \quad {\mathrm{and}}
\quad \quad \{\mu,\{N, \omega\}\} =0 \ .
$$
\begin{theorem}\label{theorem8}
If $(\omega, N)$ is an $\Omega N$-structure and $\omega$ is
non-degenerate, then $(\pi,N)$, where $\pi^\sharp = N 
\circ (\omega^\flat)^{-1}$, is a
$PN$-structure. Conversely, if $(\pi, N)$ is a $P N$-structure and 
$\pi$ is non-degenerate, then $(\omega,N)$, where $\omega^\flat =
(\pi^\sharp)^{-1} \circ N$, is an
$\Omega N$-structure.
\end{theorem}
\noindent{\it{Proof}} If $(\omega,N)$ is an $\Omega N$-structure,
we conclude from \eqref{T=0} that 
$\{[\pi,\pi]_\mu, \omega\}, \omega\} =~0$. If $\omega$ is
non-degenerate with inverse $\sigma$, applying $\{\sigma,
\cdot \}$ twice yields $[\pi,\pi]_\mu = 0$.  From \eqref{C=0} we then
see that $\pi$ and $N$ are compatible, therefore $(\pi,N)$ is a
$PN$-structure.

Conversely, if $(\pi, N)$ is a $P N$-structure, from \eqref{C=0}
we obtain $\{\{\pi,d_\mu \omega\},\pi\} = 0$.
If $\pi$ is non-degenerate with inverse $\tau$, applying 
$\{\tau,\cdot \}$ twice yields $d_\mu\omega = 0$.
Then \eqref{T=0} yields $\{\pi, d_N \omega\} = 0$. Applying
$\{\tau,\cdot \}$ yields $d_N\omega = 0$, whence also
$d(\omega_N) = 0$.
\hfill $\square$
\begin{theorem}\label{theorem9}
If $(\pi, \omega)$ is a $P\Omega$-structure, then $(\omega,N)$, where
$N = \pi^\sharp \circ \omega^\flat$, is an
$\Omega N$-structure.
Conversely, if $(\omega, N)$ is an $\Omega N$-structure and 
$\omega$ is
non-degenerate, then $(\pi, \omega)$, where 
$\pi^\sharp = N \circ (\omega^\flat)^{-1}$, 
is a $P \Omega$-structure.
\end{theorem}
\noindent{\it{Proof}}
 If $(\pi, \omega)$ is a $P\Omega$-structure,
we conclude from \eqref{T=0} that $N$ is a Nijenhuis tensor
and therefore $(\omega ,N)$ is an
$\Omega N$-structure.

Conversely, if $(\omega, N)$ is an $\Omega N$-structure, from \eqref{C=0}, 
we obtain $\{[\pi,\pi]_\mu, \omega\}, \omega\} =~0$. If $\omega$ is
non-degenerate with inverse $\sigma$, applying $\{\sigma,
\cdot \}$ twice yields $[\pi,\pi]_\mu = 0$, and therefore
$(\pi,\omega)$ is a $P \Omega$-structure.
\hfill $\square$

\medskip

We can relate the notion of $\Omega N$-structure to that of a complementary $2$-form. 
The following proposition is a consequence of Theorems
\ref{POmegacompl} 
and \ref{theorem9}.

\begin{proposition}
If $\omega$ is a $d_\mu$-closed complementary form for a Poisson
bivector $\pi$, then 
 $(\omega,N)$, where $N = \pi^\sharp \circ \omega^\flat$, is an
$\Omega N$-structure. Conversely, if  $(\omega,N)$ is an
$\Omega N$-structure and $\omega$ is non-degenerate, then 
$\pi$ such that $\pi^\sharp = N 
\circ (\omega^\flat)^{-1}$ is a Poisson bivector and $\omega$ is a
$d_\mu$-closed complementary $2$-form for $\pi$.
\end{proposition}

\begin{remark}{\rm{Magri and Morosi \cite{MM} introduced the $\Omega N$-structures as
follows:
an $\Omega N$-structure on a manifold $M$ is a
pair $(\omega, N)$, where $\omega$ is a closed $2$-form and $N$ is a
Nijenhuis tensor such that $\omega^\flat \circ N \! = \! N^* \circ
\omega^\flat$ and $S(\omega, N) \! = \! 0$,
where $S(\omega, N)$ is the section of $\wedge^2(T^*M) \otimes T^*M$ such that,
for $X$ and $Y \in {\rm{\Gamma}}(TM)$, $S(\omega, N)(X,Y) = ({\mathcal
  L}_{NY} \omega)^\flat (X) - ({\mathcal
  L}_{NX} \omega)^\flat (Y) - (\omega^\flat \circ N) [X,Y] + d
<(\omega^\flat \circ N)  (Y), X >$.
They then proved
that for any $2$-form $\omega$ and $(1,1)$-tensor $N$ on $M$,
for all vector fields $X,Y,Z$,
$$
S(\omega, N)(X,Y,Z)= d\omega(NX,Y,Z) -d\omega(X,NY,Z)
+d(\omega_N)(X,Y.Z) \ ,
$$
with $\omega_N^\flat = \omega^\flat \circ N$.
Using this formula, they proved that if $(\pi,\omega)$ is a $P \Omega$-structure on $M$,
then $(\omega,N)$, where $N= \pi^\sharp \circ \omega ^\flat$, 
is an $\Omega N$-structure.
The preceding formula remains valid
in the case of a Lie algebroid when $d$ is replaced by $d_\mu$.
Therefore Definition \ref{defomegaN} agrees with 
the original definition of Magri and Morosi.}}
\end{remark}

\section{What is a Hitchin pair on a Lie algebroid?}\label{9}

Recall that a $d_\mu$-closed non-degenerate $2$-form is called {\it
  symplectic}. The following definition is due to Crainic \cite{C}.

\begin{definition}
A symplectic form $\omega$ and a $(1,1)$-tensor $N$ define a
{\it Hitchin pair} on a Lie algebroid $(A,\mu)$
if $\omega^\flat \circ N = N^* \circ \omega^\flat$
and $d_\mu (\omega_N)=0$,
where $\omega_N^\flat = \omega^\flat \circ N$.
\end{definition}

It follows that an $\Omega N$-structure $(\omega, N)$ when $\omega$ is
non-degenerate is a Hitchin pair.
Conversely a Hitchin pair $(\omega,N)$ is an $\Omega N$-structure if and only
the Nijenhuis torsion of $N$ vanishes.

The $2$-form $\lambda = - (\omega + (\omega_N)_N)$, where $(\omega_N)_N = \frac{1}{4} i_N (i_N \omega) = \frac{1}{4} \{N,\{N,\omega\}\}$, is called the {\it twist} of the Hitchin pair.
\begin{lemma}
Let $(\omega, N)$ be a Hitchin pair.
The twist $2$-form $\lambda$ satisfies the relation, 
\begin{equation}\label{twist}
{\mathcal T}_\mu N = \{\sigma, d_\mu\lambda\} \ ,
\end{equation}
where $\sigma$ is the inverse of $\omega$.
\end{lemma}
\noindent{\it{Proof}} Since $N = \sigma^\sharp \circ
(\omega_N)^\flat$, we can apply Formula \eqref{T=0} to $\sigma$ and
$\omega_N$. Since $d_\mu \omega=0$ and therefore $[\sigma ,\sigma]_\mu
= 0$, and since $d_\mu (\omega_N) =0$, Formula \eqref{T=0} reduces to
$$
{\mathcal T}_\mu (N)  = \frac{1}{2} \{ d_N\omega_N , \sigma \} \ .
$$
Since $d_\mu (\omega_N) = 0$, $d_N (\omega_N) = - d_\mu(i_N \omega_N) = - 2 d_\mu((\omega_N)_N) = - \frac{1}{2} d_\mu(i_N (i_N \omega))$, and we obtain 
\eqref{twist} without further computations. \hfill $\square$

\medskip

In view of Lemma \ref{lemma3}, this result agrees with the computation in \cite{C} (but there is a misprint in the definition of the twist in the unpublished preprint \cite{C}).

When $(\omega, N)$ is a Hitchin pair on a Lie algebroid $A$ 
satisfying the algebraic conditions, $N^2 - \sigma^\sharp \circ
\lambda^\flat = - {\rm{Id}}_{A}$, where $\lambda$ is the twist
$2$-form, then 
${\mathcal N} = \sigma  + N + \lambda$ is a generalized complex structure (see Section \ref{courant}) on $A \oplus A^*$, {\it{i.e.}}, 
on $M$ if $A =TM$. In matrix form,
$${\mathcal N} = 
\begin{pmatrix}
N & \sigma  \\
\lambda & - N^*
\end{pmatrix}  \ .
$$
Then $(\sigma ,N, - d_\mu \lambda)$ is a Poisson quasi-Nijenhuis structure \cite{SX} \cite{A}. 
In this case, we obtain an alternate proof of \eqref{twist}
\bigskip

The table in the next section summarizes the main definitions and implications of the preceding sections.

\bigskip
In the diagram that summarizes the relationships between $PN$- ,
$P\Omega$- and $\Omega N$-structures, the arrows denote implications, and the
dotted arrows denote implications under a non-degeneracy assumption.

\newpage

\section{Summary}\label{table}

\subsection{Definitions}
$PN$ ($NP=PN^*$)
$$
\boxed{ \{\{\pi,\mu\},\pi\} = 0, \quad \{\{\pi,\mu\},N \} + \{\{N,\mu\},\pi\}
= 0, \quad \{N, \{N,\mu\}\} - \{ N^2,\mu\}=0}
$$

\noindent
$P\Omega$                 
$$
\boxed{ \{\{\pi,\mu\},\pi\} = 0, \quad \{\mu, \omega\} =0, \quad  \{\{\{\pi,\omega\},
  \mu\}, \omega\} =0}
$$

\noindent
$\Omega N$ ($\omega^\flat \circ N = N^* \circ \omega^\flat$)
$$
\boxed{ \{\mu, \omega\} =0, \quad  \{N, \{N,\mu\}\} - \{ N^2,\mu\} = 0, \quad 
  \{\mu, \{N,\omega\}\}= 0}
$$

\noindent
Hitchin pair ($\omega^\flat \circ N = N^* \circ \omega^\flat$)
$$
\boxed{ \{\mu, \omega\} =0, \quad 
  \{\mu, \{N,\omega\}\}= 0}
$$

\noindent
Complementary $2$-form
$$
\boxed{ \{\{\pi,\mu\},\pi\} = 0, \quad  \{\{\omega,\{\pi,\mu\}\},\omega\} = 0 }
$$

\medskip

\subsection{Relationships}

$$
\boxed{P\Omega \Longrightarrow PN \quad (N = \pi \circ \omega)}
$$
$$          
\boxed{PN \quad and \quad \pi \quad non{\mathrm{-}}degenerate \Longrightarrow 
P\Omega \quad (\omega = \pi^{-1} \circ N)}
$$
$$
\boxed
{\Omega N \quad and \quad \omega \quad non{\mathrm{-}}degenerate \Longrightarrow
  PN 
\quad (\pi = N \circ \omega^{-1})}
$$
$$          
\boxed{PN \quad and \quad \pi \quad non{\mathrm{-}}degenerate \Longrightarrow 
\Omega N \quad (\omega = \pi^{-1} \circ N)}
$$
$$
\boxed
{P\Omega \Longrightarrow \Omega N \quad (N=\pi \circ \omega)}
$$
$$
\boxed
{\Omega N \quad and \quad \omega \quad non{\mathrm{-}}degenerate \Longrightarrow
  P\Omega \quad 
(\pi = N \circ \omega^{-1})}
$$
$$
\boxed
{Hitchin \quad pair \quad and \quad N \quad Nijenhuis \Longleftrightarrow
\Omega N \quad and \quad \omega \quad non{\mathrm{-}}degenerate }
$$
$$
\boxed
{\omega \quad closed \quad complementary \quad 2{\mathrm{-}}form \quad
  for \quad \pi \Longleftrightarrow P\Omega}
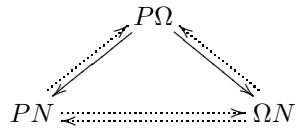
$$

\hspace*{3,5cm}
 \xymatrix{
   & P\Omega \ar[ld] \ar[rd]&  \\
 PN \ar@{.>}[ur]<3pt> \ar@{.>}[rr] &  & \Omega N \ar@{.>}[ul]<-3pt> \ar@{.>}[ll]<3pt>
  }

\newpage

\section{Nijenhuis tensors on Courant algebroids}\label{courant}

\subsection{Courant algebroids}
The double $E = A \oplus A^*$ of any proto-bialgebroid $(A,A^*)$ is a Courant
algebroid when equipped with the Dorfman bracket (see Section \ref{sectiondorfman}) 
which is the derived
bracket defined by
$$
[u,v]_S = \{\{u,S\},v\} \ ,
$$
for $u, v \in {\rm{\Gamma}} (A \oplus A^*)$. Here $S$ denotes the structure
$\phi + \mu + \gamma + \psi$. When $(A,A^*)$ is a Lie bialgebroid, then
$S = \mu + \gamma$ and $[~,~]_S$ is the bracket that we have denoted
by $[~,~]_D$ in Section \ref{sectiondorfman}. The usual case is $A = TM$ with $\gamma =0$.
(See \cite{Gua}.)

The skew-symmetrization of the Dorfman bracket is the Courant bracket,
$[u,v]'_S = \frac{1}{2} ([u,v]_S - [v,u]_S)$.

\subsection{Nijenhuis tensors and generalized complex structures}
For an endomorphism $\mathcal N$ of $E$, the Dorfman torsion ${\mathcal
  T}_S{\mathcal N}$ (resp., the Courant torsion ${\mathcal
  T}'_S{\mathcal N}$) of ${\mathcal N}$ is defined by
$$
({\mathcal T}_S{\mathcal N})(u,v)=
[{\mathcal N}u,{\mathcal N}v]_S - {\mathcal N}([{\mathcal N}u,v]_S
+[u,{\mathcal N}v]_S)
+{\mathcal N}^2[u,v]_S
$$
(resp., $
({\mathcal T}'_S{\mathcal N})(u,v)=
[{\mathcal N}u,{\mathcal N}v]'_S - {\mathcal N}([{\mathcal N}u,v]'_S
+[u,{\mathcal N}v]'_S)
+{\mathcal N}^2[u,v]'_S$).

A \emph{generalized almost complex structure} is an endomorphism ${\mathcal
  N}$ of $E$ which is orthogonal
with respect to the symmetric bilinear form, and
such that ${\mathcal N}^2= - {\mathrm{Id}}_E$.
Note that if ${\mathcal N}^2 = - \lambda {\rm{Id}}_E$, with $\lambda
\neq 0$, then
${\mathcal N}{\mathcal N}^*= {\rm{Id}}_E$ is equivalent to ${\mathcal N} +
\lambda 
{\mathcal N}^* = 0$. In particular, the orthogonality condition for 
a generalized
almost complex structure is equivalent to ${\mathcal N} + {\mathcal N}^*=0$.

A generalized almost complex structure ${\mathcal N}$
is a \emph{generalized complex structure} if
${\mathcal T}_S{\mathcal N}= 0$.

\begin{lemma} For a generalized almost complex structure ${\mathcal
    N}$, ${\mathcal T}_S{\mathcal N}= {\mathcal T}'_S{\mathcal N}$.
\end{lemma}

\noindent{\it Proof}
Since ${\mathcal N}$
is orthogonal, the equality follows from the relation
$[u,v]_S + [v,u]_S
= {\mathcal D}(<u,v>)$, where $({\mathcal D} u) (f)= 
<\rho u, df>$. \hfill $\square$

\medskip

Thus the integrability condition for an almost complex structure can be expressed either in terms of non-skew symmetric
brackets or in terms of skew-symmetric brackets. 

As observed by Grabowski in \cite{G}, Formula \eqref{a} is
valid when $\mu$ denotes the cubic homological function defining the Courant
algebroid structure, denoted by $S$ above.
Therefore the following theorem is proved.

\begin{theorem}\label{GCS}
The generalized complex structures are the generalized almost
complex structures such that
$$
\{\{{\mathcal N}, S \}, {\mathcal N}\} = S \ .
$$
\end{theorem}

The preceding remarks also apply to almost product and almost subtangent structures. They are characterized by 
$\{\{{\mathcal N}, S \}, {\mathcal N}\} = - S$
and
$\{\{{\mathcal N}, S \}, {\mathcal N}\} = 0$, respectively.

When the torsion of ${\mathcal N}$ vanishes, ${\mathcal N}$ defines a
{\it{deformed structure}} $S_{\mathcal N} =\{{\mathcal N}, S\}$
on~$A$, and ${\mathcal N}$ maps bracket $[~,~]_{S_{\mathcal N}}$ to 
bracket $[~,~]_{S}$. A necessary and sufficient condition for 
$S_{\mathcal N} =\{{\mathcal N}, S\}$ to be a Dorfman bracket is 
$\{S,\{\{{\mathcal N}, S\}, {\mathcal N} \}\} = 0$. In fact, Formula \eqref{b}
extends to the Courant algebroid case, 
$$
\frac{1}{2}\{S_{{\mathcal N}},
S_{{\mathcal N}}\}
= \frac{1}{2}\{\{{\mathcal N},S\},\{{\mathcal N},S\}\}
=\{S,{\mathcal T}_S{\mathcal N}\} \ .
$$
The vanishing of ${\mathcal T}_S({\mathcal N})$ expresses the fact that 
$ {\mathcal N}:(E, S_{{\mathcal N}}) \to (E , S)$ preserves the brackets.

\medskip

In the next sections we shall study Monge-Amp\`ere structures as examples of
compatible structures and generalized geometries.

\section{Monge-Amp\`ere structures on manifolds}\label{MA}

The notion of Monge-Amp\`ere
structure has its origin in the theory of symplectic Monge-Amp\`ere
operators and equations. 
See \cite{KLR} for a detailed analysis of symplectic and contact Monge-Amp\`ere operators and equations, together with many examples.
Let $M$ be a smooth manifold of dimension $n$ and let $T^*M$ be its
cotangent bundle. 
We shall denote the space of $k$-forms 
on $T^*M$ by ${\rm{\Omega}}^k(T^*M)$
and the space of vector fields by
${\mathcal X}(T^*M)$. Let ${\mathcal F}(T^*M)$ be the space of functions on
the supermanifold $T^*[2](T(T^*M))[1]$.

We shall denote the canonical symplectic 2-form on
$T^*M$ by $\Omega$,
and its inverse, the canonical bivector, by $\pi_\Omega$.
More generally, we shall denote by $\pi_\tau$ the bivector on $T^*M$
that is the 
inverse of a non-degenerate 2-form $\tau$ on $T^*M$.

\begin{definition}
The pair $(\Omega,\omega)$ is a {\emph {Monge-Amp\`ere
    structure}} on $M$ if $\omega$ is an $n$-form on $T^*M$ satisfying the condition
$\omega\wedge \Omega =0.$
\end{definition}

 According to the original ideas of Lychagin
\cite{Ly} \cite{KLR}, any symplectic Monge-Amp\`ere operator on $M$
can be defined by an
{\emph{effective form}} on $T^*M$ of degree $k$, $2 \leq k\leq n$,
{\it i.e.}, a $k$-form $\omega$ such that  $i_{\pi_{\Omega}}\omega = 0$. 
%A  $k$-form is primitive in this sense if and only if it is
%primitive in the sense of Definition \ref{prim} of Section
%\ref{representation}.
When $k=n$, this condition is equivalent to
the condition $\omega\wedge\Omega = 0$, so $(\Omega, \omega)$ is a
Monge-Amp\`ere structure if and only if $\omega$ is an effective
$n$-form.

A correspondence between forms on $T^*M$ and form-valued differential
operators on $M$
is defined as follows.
Let $k$ be a positive integer, $k \leq n$.
For any $k$-form  $\omega$ on $T^*M$, define the {\it {Monge-Amp\`ere
    operator}}, ${\rm{\Delta}}_{\omega} : C^{\infty}(M)\to
{\rm{\Omega}}^k(T^*M)$,
by
\begin{equation}\label{MAO}
{\rm{\Delta}}_{\omega}(f) = (df)^*(\omega) \ ,
\end{equation}
for any $f\in C^{\infty}(M)$.
We understand the differential $df$ to be a map, $df : M\to T^*M$, the
section of the cotangent 
bundle defined by the smooth function $f$. 
The equation ${\rm{\Delta}}_{\omega}(f) = 0$ is called a {\it {Monge-Amp\`ere equation}}.

\section{Monge-Amp\`ere structures in dimension $2$ and compatible structures}\label{MA2}
\subsection{Monge-Amp\`ere structures in 
dimension $2$}

We first consider the simplest geometric examples, those of
Monge-Amp\`ere structures in dimension $n=2$.
In this case, $\dim (T^*M )= 4$, and a Monge-Amp\`ere structure is
defined by a pair of $2$-forms $(\Omega,\omega)$ on $T^*M$,
where $\Omega$ is the canonical $2$-form and $\omega$ satisfies the
effectivity condition, $\omega\wedge \Omega = 0$, or equivalently,
$i_{\pi_{\Omega}}\omega =0$.
We consider the $(1,1)$-tensor on $T^*M$, $A_{\omega} =
(\pi_\Omega)^\sharp \circ \omega^\flat$, which satisfies, for all $X,Y \in {\mathcal X}(T^*M)$,
$$
\omega(X,Y) = \Omega(A_{\omega}X,Y) \ .
$$
It is easy to prove that $A_{\omega}$ satisfies the equation
$A_{\omega}^2 + \Pf(\omega){\rm Id}=0$, where
the Pfaffian, $\Pf(\omega)$, of the 2-form $\omega$ is defined by
$$
\Pf(\omega)\Omega\wedge\Omega = \omega\wedge\omega \ ,
$$
and ${\rm Id}$ is the identity of $T(T^*M)$.

\subsection{Non-degenerate Monge-Amp\`ere structures in
dimension $2$}
By definition, a Monge-Amp\`ere structure $(\Omega, \omega)$ on 
$T^*({\mathbb R}^2)$ is called
{\it non-degenerate} if its Pfaffian $\Pf(\omega)$ is nowhere-vanishing.

When $(\Omega, \omega)$ is a non-degenerate Monge-Amp\`ere structure,
we consider
the normalized 2-form ${\widetilde{\omega}} =
\frac{\omega}{\sqrt{|\Pf(\omega)|}}$, with inverse bivector
$\pi_{{\widetilde{\omega}}}=
\sqrt{|\Pf(\omega)|}\pi_\omega$. The normalized $(1,1)$-tensor
$J_{\omega}$ is defined by
$J_{\omega} =
(\pi_\Omega)^\sharp \circ {\widetilde \omega}^\flat$, and it satifies
$$
J_{\omega} = \frac{A_{\omega}}{\sqrt{|\Pf(\omega)|}} \ .
$$
Then $J_{\omega}^2 =-  {\rm Id}$ or $J_{\omega}^2 =  {\rm Id}$.
The sign of $\Pf(\omega)$ determines whether
the corresponding Monge-Amp\`ere
  operator is elliptic (when $\Pf(\omega) > 0$ and therefore 
$J_\omega^2 = -  {\rm Id}$) or hyperbolic (when $\Pf(\omega) < 0$ and therefore
  $J_\omega^2 = {\rm Id})$.

It is proved in \cite{KLR} that the integrability of $J_{\omega}$,
{\it{i.e.}} the condition ${\mathcal
  T}_\mu({J_{\omega}})=0$, where $\mu\in {\mathcal F}^{0,1}(T^*M)$
defines the standard Lie algebroid structure of $T(T^*M)$,
is equivalent to the
condition that the normalized 2-form ${\widetilde{\omega}}$ be closed.
This integrability condition is also equivalent to the existence of a
symplectomorphism mapping the $2$-form $\omega$ to a form with
constant coefficients. The corresponding Monge-Amp\`ere
operator ${\rm{\Delta}}_{\omega}$ is then equivalent to an
operator with constant coefficients.

\subsection{Properties of non-degenerate Monge-Amp\`ere 
structures in dimension $2$}\label{MA2modular}
We show that in the case of dimension $2$,
non-degenerate Monge-Amp\`ere\ structures give rise to composite
structures.

{$\bullet$} By Theorem \ref{converse} of Section \ref{sectioncompatible},
if $(\Omega,\omega)$ is a
non-degenerate Monge-Amp\`ere\ structure on $M$ satisfying the
condition
$d{\widetilde{\omega}} = 0$, then
the pairs $(\pi_{\Omega}, J_{\omega})$ and
$(\pi_{\widetilde \omega}, J_{\omega})$ are
$PN$-structures on $T^*M$, {\it{i.e.}}, on the Lie algebroid $T(T^*M)$.

$\bullet$
Theorem \ref{theorempomega}
in Section \ref{sectionpomega} implies that, when $d{\widetilde \omega}
=0$, the pair
$(\pi_{\Omega},\widetilde{\omega})$ is a $P\Omega$-structure
and the pair $(\widetilde{\omega}, J_{\omega})$ 
is an $\Omega N$-structure on $T^*M$.

$\bullet$ 
  Let $\mu_{J_\omega}$ be the element of ${\mathcal F}^{0,1}(T^*M)$ defined by
$$
\mu_{J_\omega} = \{J_{\omega},\mu\} \ ,
$$
where as above $\mu$ is the Lie algebroid structure of $T(T^*M)$.
When $J_{\omega}$ is integrable,
$\mu_{J_\omega}$ defines
a new Lie algebroid structure on $T(T^*M)$ deformed by $J_\omega$. By
Remark \ref{remarkcompatible} of Section \ref{sectionNij}, this
structure is compatible with the standard structure,
$$
\{\mu + \mu_{J_\omega},\mu + \mu_{J_\omega}\} =0 \ .
$$

The observation in Section \ref{modclass} can be applied to the
modular class of this deformed Lie algebroid structure. Assume that
$\widetilde \omega$ is closed. Then the deformed structure 
$\mu_{J_{\omega}}$ is equal to $\{\{\pi_{\Omega},
\mu\},\widetilde \omega\}$. In fact, by the Jacobi identity, since 
$\{\mu,\widetilde \omega\}=0$,
$$
\mu_{J_\omega} =\{J_{\omega},\mu\}
= \{\{\pi_{\Omega},
\widetilde \omega\},\mu\}
= \{\{\pi_{\Omega},\mu\},\widetilde \omega\} \ .
$$
Hence, the $1$-form on $T^*M$, $\xi_{\pi_{\Omega},\widetilde\omega,\lambda},$
defined by the Liouville volume form $\lambda = \frac{1}{2}\Omega \wedge \Omega$, satisfying
$$
\partial^{\pi_\Omega}_{\widetilde\omega,\lambda}
- \partial^{\pi_\Omega}_{\widetilde\omega} =
i_{\xi_{\pi_{\Omega},\widetilde\omega,\lambda}}
$$
is a $d_{J_\omega}$-cocycle whose cohomology class is the modular class of
the Lie algebroid $(T(T^*M),\mu_{J_\omega})$ defined by the
Nijenhuis operator $J_{\omega}$.
\begin{proposition}\label{unimod}
The deformed Lie algebroid $(T(T^*M),\mu_{J_\omega})$ which is obtained from
a Monge-Amp\`ere structure such that $d\widetilde \omega =0$
is unimodular.
\end{proposition}
\noindent{\it Proof}
Since ${J_\omega}$
is a Nijenhuis tensor, 
the modular class in the $d_{J_\omega}$-cohomology of the Lie 
algebroid $(T(T^*M),\mu_{J_\omega})$ is the class of the $1$-form $d({\rm Tr}J_\omega)$ (see \cite{KM1} \cite{DF}).
Since the form $\omega$ is effective, the $(1,1)$-tensor
$A_\omega$, and hence $J_\omega$, are traceless. \hfill $\square$

\medskip

$\bullet$ 
Any non-degenerate
Monge-Amp\`ere structure
  $(\Omega,\omega)$ on $M$ such that $\omega$ is closed defines a Hitchin pair
  $(\Omega,A_{\omega})$ in the sense of Crainic \cite{C} on $T^*M$. If, in particular, 
  this
  structure is defined by a non-degenerate Monge-Amp\`ere operator
  with constant coefficients,
  the $(1,1)$-tensor $A_{\omega}$ is integrable and
  $(\Omega,A_{\omega})$ is an
$\Omega N$-structure on $T^*M$.

$\bullet$ Monge-Amp\`ere structures of divergence type were
  defined in \cite{KLR}. A pair $(\Omega, \omega)$, where $\omega$ is
  a $2$-form, is called a structure of divergence type if there exists 
  a function $\phi$
on $T^*M$ such that $\omega + \phi  \Omega$ is closed. 
  Following \cite{B}, we observe that a non-degenerate
 structure $(\Omega,\omega)$ of divergence type, where $\omega$ is not
  necessarily effective, defines a generalized almost complex structure
${\mathcal J}_{\omega} =
    \begin{pmatrix}
    A_{\omega}& \pi_{\Omega}^\sharp\\
    -\Omega^\flat ({\mathrm{Id}}+A^2_{\omega})& -A^*_{\omega}\\
    \end{pmatrix}$ on $T^*M$.  The pair $(\Omega, A_\omega)$ is a Hitchin pair if and only 
    if ${\mathcal J}_\omega$ is integrable if and only if $\omega$ is closed.  
    
If, in addition, the Monge-Amp\`ere structure
$(\Omega,\omega)$ satisfies the condition $d\widetilde\omega = 0$, 
    we obtain another generalized complex structures on $T^*M$,
   ${\mathbb J}_{\omega} =
    \begin{pmatrix}
    J_{\omega}& \pi_{\Omega}^\sharp\\
    0& -J^*_{\omega}\\
    \end{pmatrix}$
if ${\rm \Delta}_{\omega}$ is equivalent to an elliptic Monge-Amp\`ere operator
    with constant coefficients, and
    ${\mathbb J}'_{\omega} =
    \begin{pmatrix}
    J_{\omega}& \pi_{\Omega}^\sharp\\
    -2\Omega^\flat& -J^*_{\omega}\\
    \end{pmatrix}$
which corresponds to a hyperbolic Monge-Amp\`ere operator.

The tensor ${\mathbb J}_{\omega}$ can be written as ${\mathbb J}_{\omega}
= \pi_\Omega + J_{\omega}$ since
$$
{\mathbb J}_{\omega} (u) = \{u, \pi_\Omega + J_{\omega}\} \ ,
$$
for all $u \in {\mathcal F}(T^*(M)$, and similarly for ${\mathcal
  J}_\omega$ and ${\mathbb J'}_\omega$.

$\bullet$ The deformed Lie bialgebroid structure on $(T(T^*M),T^*(T^*M))$ defined by
  ${\mathbb J}_{\omega}$ induces a new Courant algebroid structure on $T(T^*M)
  \oplus T^*(T^*M)$ which we shall call a
{\it{Monge-Amp\`ere Courant algebroid}}.
This structure is defined by the function $S_{{\mathbb J}_{\omega}} =
\{{\mathbb J}_{\omega},S\}
\in {\mathcal F}(T^*M)$, 
where $S = \mu$ is the standard Courant algebroid structure of $T(T^*M)
  \oplus T^*(T^*M)$.

The integrability condition of Theorem \ref{GCS} is
$$
\{\{{\mathbb J}_{\omega},S\},{\mathbb J}_{\omega}\} = S \ .
$$
When the integrability condition is satisfied, $S_{{\mathbb J}_{\omega}}$
satisfies 
$\{S_{{\mathbb J}_{\omega}},
S_{{\mathbb J}_{\omega}}\}=0$ and $S_{{\mathbb J}_{\omega}}$
maps the Dorfman bracket defined by $S_{{\mathbb J}_{\omega}}$
to the bracket defined by $S$.

\subsection{The von Karman equation}
The conditions $d\omega = 0$ and
$d(\frac{\omega}{\sqrt{|\Pf(\omega)|}})=0$ are very different. In the
former case, there is a pair of symplectic forms on $T^*M$. The
latter condition is the necessary and sufficient condition for the
Monge-Amp\`ere structure to be equivalent to a
structure with constant coefficients.
If the Monge-Amp\`ere structure $(\Omega, \omega)$ is
  such that  $d{\widetilde{\omega}}  \neq 0$, then the
  torsion of $J_{\omega}$ does not vanish, the integrability condition
  is not satisfied.
The following example shows
that the condition
$d\omega=0$ is not sufficient to define a $PN$- or an $\Omega N$-structure.

Let $(q_1,q_2,p_1,p_2)$ be the canonical coordinates on $T^*({\mathbb
  R}^2) ={\mathbb R}^4$.
Let $(\Omega,\omega)$ be the Monge-Amp\`ere structure on ${\mathbb R}^2$ defined by the $2$-form on $T^*({\mathbb R}^2)$,
$$
\omega = p_1 dp_1 \wedge dq_2  - dp_2 \wedge dq_1 \ .
$$
The corresponding partial differential equation is the  {\it von Karman equation},
$$
f_{q_1}f_{q_1 q_1} - f_{q_2 q_2} = 0 \ .
$$
 It is easy to show that $\omega \wedge \Omega = 0$ and $d\omega =
 0$.
 This structure is
 non-degenerate
 in the complement of the hyperplane $p_1 = 0$, since the Pfaffian $\rm
 Pf(\omega)$ is equal to $p_1$.
 In the half-space $p_1 > 0$ (resp.,  $p_1 < 0$)
the von Karman equation is an elliptic (resp., hyperbolic)
Monge-Amp\`ere equation. Now  $d{\widetilde{\omega}} \neq 0$ since
\begin{equation}
  d{\widetilde{\omega}} = d(\frac{\omega}{\sqrt{|\rm Pf (\omega)|}})=
  \frac{1}{2} |p_1|^{-3/2} dp_1 \wedge dp_2 \wedge dq_1 \ .
\end{equation}
Therefore,
the Monge-Amp\`ere structure $(\Omega,\omega)$ is not equivalent
 to a Monge-Amp\`ere structure with constant coefficients.
 It  does not define a $PN$- nor an $\Omega N$-structure on
${\mathbb R}^4$, nor a Monge-Amp\`ere Courant
 structure because the equation
 $\{S_{{\mathbb J}_{\omega}}, S_{{\mathbb J}_{\omega}} \}=0$ is not
 satisfied.
The Poisson
 tensor inverse to $\omega$ is
 $$
 \pi_{\omega}=\frac{1}{|p_1|}\frac{\partial}{\partial p_1}\wedge
 \frac{\partial}{\partial q_2} - \frac{\partial}{\partial p_2}\wedge
 \frac{\partial}{\partial q_1} \ .
 $$
The computation of
$\{\{\pi_{\Omega},\mu\},\pi_{\omega}\}$ shows that
the Schouten bracket of $\pi_{\omega}$ and the canonical Poisson
tensor $\pi_{\Omega}$ is the $3$-vector,
 $$
 [\pi_{\Omega},\pi_{\omega}]_\mu = -\frac{1}{(p_1)^{2}}\frac{\partial}{\partial
   q_1}\wedge\frac{\partial}{\partial
   q_2}\wedge\frac{\partial}{\partial p_1} \ .
 $$

\subsection{Generalized
 Monge-Amp\`ere structures}

More generally, we can consider \emph{generalized
 Monge-Amp\`ere structures} $(\omega_1,\omega_2)$, where both $2$-forms
on $T^*M$, $\omega_1$ and $\omega_2$,
 are non-degenerate but not necessarily closed. 
 The corresponding equations are
 systems of non-linear first-order  partial differential equations
whose non-linearity has a
specific form. Such systems, called \emph{Jacobi systems}, are studied in
\cite{KLR}. 
A Jacobi system is called {\it{non-degenerate}} if $\omega_1 \wedge \omega_2 = 0$ 
and there exists  a nowhere vanishing function on $T^*M$, 
$\epsilon$, such that $\omega_1 \wedge \omega_1=\epsilon \, \omega_2 \wedge \omega_2$.
The Jacobi systems are of the form, for a pair of functions $(u,v)$ on $M={\mathbb R}^2$
with coordinates $(x, y)$,
\begin{equation}\label{jacob}
\begin{cases} 
a+b u_x +c u_y +
d v_x + e v_y +
f {\mathcal J}_{u,v} = 0 \ , \\
A+B u_x +C u_y +
D v_x + E v_y +
F {\mathcal J}_{u,v} = 0 \ ,\\
\end{cases}
\end{equation}
where
${\mathcal J}_{u,v}$ is the Jacobian determinant of $(u,v)$.

The Jacobi systems can be defined invariantly as follows. Let ${\cal M} = M\times {\mathbb R}^2$,
where $M$ is a 2-dimensional manifold and
let $\omega_{i}$, $i=1,2$, be 2-forms on $\cal M$. We define 
the differential operators,
${\rm \Delta}_{\omega_i} : C^{\infty}(M, {\mathbb R}^2) \to {\rm{\Omega}}^2(T^*M)$, by 
\begin{equation}\label{jacop}
{\rm \Delta}_{\omega_i}(f) = \omega_i |_{L_f} \ , \quad i=1,2 \ ,
\end{equation}
where $L_f$ is the graph of the ${\mathbb R}^2$-valued function $f$ on $M$, a $2$-dimensional
surface in $\cal M$. The system \eqref{jacob} is then written 
$${\rm \Delta}_{\omega_i}f = 0,\quad i=1,2,\quad f=(u,v).$$
If the restrictions of $\omega_1$ and $\omega_2$ to
the surface $L_f \subset {\cal M}$ vanish, we shall say that $L_f$ is a {\it{generalized solution}} of \eqref{jacob}.
Geometrically, we can assign to each point $m\in {\cal M}$ the plane
in $\wedge^2(T^*_m{\mathcal M})$, called the {\it Jacobi plane}, generated by
$\omega_1|_m$ and $\omega_2|_m$, thus defining a smooth distribution
on $\mathcal M$
which corresponds to the system \eqref{jacop}. The submanifold $L_f$ 
is an integral manifold for this distribution.

We define a $(1,1)$-tensor $A\in {\rm \Gamma}(T{\cal M} \otimes T^*{\cal M})$ by
$$
\omega_2(X,Y) = \omega_1(AX,Y) \ ,
$$
for all $X$ and $Y \in {\rm \Gamma}(T{\cal M})$.
If the Jacobi system \eqref{jacop} is non-degenerate and if, in
addition, $\epsilon= 1$ or $\epsilon= - 1$,
then $A^2 = \epsilon$ and we can associate to such systems an almost complex or almost product structure on $T\cal M$
(see \cite{KLR}).
Let $\pi_{\omega_i}\in 
\wedge^2(T\cal M)$, $i=1,2$, be the bivectors which are inverse to the
non-degenerate $2$-forms $\omega_i$.
Suppose that these bivectors satisfy the following conditions,
\begin{equation}{\label{bihamilt}}
\begin{cases}
[\pi_{\omega_1},\pi_{\omega_1}]=[\pi_{\omega_2},\pi_{\omega_2}] \ ,&\\
[\pi_{\omega_1},\pi_{\omega_2}]=0 \ .&\\
\end{cases}
\end{equation}
A pair of bivectors satisfying the conditions \eqref{bihamilt} is called a
{\it Hitchin pair of bivectors} in \cite{B}.
Theorem \ref{Nij} of Section \ref{sectioncompatible}
implies that the Jacobi systems on ${\mathbb R}^2$
associated with Hitchin
pairs of Poisson bivectors define $PN$-structures on ${\mathbb R}^4$.

\section{Monge-Amp\`ere structures in dimension $3$ and generalized
  geometry}\label{MA3}

\subsection{Classification}
Dimension $3$ plays an exceptional role in the geometry of
Monge-Amp\`ere operators. The classification problem for
Monge-Amp\`ere operators and equations
on $3$-dimensional manifolds
can be reduced to a classical problem in geometric invariant theory, the
determination of the
normal forms of the effective $3$-forms in a $6$-dimensional real
symplectic vector space $V$, in other words, of the orbits of the
symplectic group ${\rm{Sp}}(6)$ on the space
of effective $3$-forms on V.
This problem was solved in \cite{LyR} (see also
\cite{KLR}).
There are three types of generic orbits, each with a non-trivial
stabilizer, each corresponding to a non-degenerate Monge-Amp\`ere
structure with a
non-degenerate non-linear
Monge-Amp\`ere operator.
Let $(q_1,q_2,q_3,p_1,p_2,p_3)$ be the canonical coordinates on $T^*({\mathbb
  R}^3) ={\mathbb R}^6$, and let $u$ be a function on $T^*({\mathbb R}^3)$.
The three types of generic orbits are those of the following $3$-forms 
with constant coefficients, with corresponding Monge-Amp\`ere equations: 
\begin{equation}\label{hess}
\omega = dp_1\wedge dp_2\wedge dp_3 - dq_1\wedge dq_2\wedge dq_3 \ ,
 \quad \quad {\rm{\Delta}}_{\omega} = \hess(u) -1 \ ,
\end{equation}
\begin{equation}\label{SL}
\omega = dp_1\wedge dp_2\wedge dp_3 - dp_1\wedge dq_2\wedge dq_3 -
dq_1\wedge dp_2\wedge dq_3 - dq_1\wedge dq_2\wedge dp_3 \ ,
{\rm{\Delta}}_{\omega} = \hess(u) -{\rm{\Delta}}(u) \ ,
\end{equation}
\begin{equation}\label{psSL}
\omega = dp_1\wedge dp_2\wedge dp_3 - dp_1\wedge dq_2\wedge dq_3 -
dp_2\wedge dq_1\wedge dq_3- dp_3\wedge dq_1\wedge dq_2 \ ,
{\rm{\Delta}}_{\omega} = \hess(u) -\Box(u) \ ,
\end{equation}
where ${\rm{\Delta}}= \frac{\partial^2 }{\partial q_1
  ^2}+\frac{\partial^2 }{\partial q_2 ^2}+\frac{\partial^2 }{\partial
  q_3 ^2}$ is the Laplacian, $\Box = \frac{\partial^2 }{\partial q_1
  ^2}+\frac{\partial^2 }{\partial q_2 ^2}-\frac{\partial^2 }{\partial
  q_3 ^2}$ is the D'Alembertian of signature $(2,1)$,
and $\hess(u) = {\rm det}(u_{q_i q_j}), 
1\leq i,j \leq 3$, is the Hessian of the function $u$, {\it i.e.},
the determinant of the matrix of second-order partial derivatives of
$u$ with respect to $q_1,q_2,q_3$.

We shall show that, in full analogy to the $2$-dimensional case
where almost complex (resp., almost product) 
structures\footnote{In \cite{BTh}, Banos called these structures
``generalized Calabi-Yau structures'', but this terminology conflicts
with Hitchin's in
\cite{H}. Below we shall clarify the difference between these two
generalizations of the Calabi-Yau structures.} correspond
to elliptic (resp., hyperbolic) Monge-Amp\`ere operators,
there exist three generalized structures in the sense of Grabowski
\cite{G} corresponding to the three types of Monge-Amp\`ere structures
in dimension $3$.

\subsection{Hitchin endomorphism and Hitchin Pfaffian}
To each Monge-Amp\`ere structure $(\Omega,\omega) \in
{\rm{\Omega}}^2(T^*M) \times {\rm{\Omega}}^3(T^*M)$ on a
$3$-dimensional manifold $M$, where $\omega$ is effective,
are associated the following \cite{H0} \cite{KLR}:

$\bullet$
the Liouville volume form, ${\rm vol}$,
associated with $\Omega$,
$$
{\rm vol}=-\frac{1}{6} \Omega \wedge \Omega \wedge \Omega  \in
{\rm{\Omega}}^6(T^*M) \ ,
$$

$\bullet$
the Hitchin endomorphism, $H_{\omega} : {\mathcal X}(T^*M) \to {\mathcal X}(T^*M)$,
defined by
$$
H_{\omega}(X) = i_{X}\omega\wedge\omega\in {\rm{\Omega}}^5(T^*M) \simeq
 {\mathcal X}(T^*M) \ ,
$$
for all $X \in {\mathcal X}(T^*M)$, where ${\rm{\Omega}}^5(T^*M)$ is identified with
${\mathcal X}(T^*M)$ by means of the Liouville form,

$\bullet$
the Hitchin Pfaffian, $\lambda_{\omega}$, defined by
$$
\lambda_{\omega} = \frac{1}{6}{\rm Tr}(H_{\omega}^2) \ ,
$$

$\bullet$
the symmetric bilinear form, $q_\omega$, defined by
$$
q_{\omega}(X,Y) = \Omega(H_{\omega}X,Y) \ ,
$$
for all $X$ and $Y \in  {\mathcal X}(T^*M)$.

The Hitchin endomorphism and the Hitchin Pfaffian are related by
$$
H_{\omega}^2  = \lambda_{\omega} \rm Id.
$$

By definition, a Monge-Amp\`ere structure $(\Omega, \omega)$ on $T^*({\mathbb R}^3)$ is called
{\it non-degenerate} if its Hitchin Pfaffian 
$\lambda_{\omega}$ is nowhere-vanishing.

An essential part of the proof of the above-mentioned classification
is the proof that the forms in the orbit of the form $\omega$ of (\ref{hess}) have
negative Hitchin Pfaffian, while those in the orbits of the forms of (\ref{SL}) and of
(\ref{psSL}) have positive Hitchin Pfaffian
and quadratic forms $q_{\omega}$ of different signatures.

For a $2$-form $\tau$, we define the modified Pfaffian ${\mathcal{P}}f(\tau)$
by
$$
\tau \wedge \tau \wedge \Omega = - \frac{1}{3} {\mathcal{P}}f(\tau)
\Omega \wedge \Omega \wedge \Omega \ .
$$ 
The following statement is the result of a straightforward computation.
\begin{proposition}\label{LRInv}
The modified Pfaffian ${\mathcal{P}}f$, the Hitchin endomorphism $H_\omega$ and 
the bilinear form $q_\omega$ satisfy the relations
$$
{\mathcal{P}}f(i_X \omega) = \Omega(H_{\omega}X,X) = q_{\omega}(X,X) \equiv -\frac{1}{4}\imath_{\pi_{\Omega}}\imath_{\pi_{\Omega}}(i_X \omega\wedge i_X\omega) \ ,
$$
for all $X \in {\mathcal X}(T^*M)$.
\end{proposition}

\subsection{Properties of non-degenerate Monge-Amp\`ere structures in
dimension $3$}
We can now draw conclusions analogous to those of the  $2$-dimensional
case of Section~\ref{MA2}.

$\bullet$ Any Monge-Amp\`ere structure $(\Omega,\omega)$
satisfies the conditions of
Lemma \ref{lemma3} of Section \ref{bivectors} and therefore,
$$
\{\{ X , \{ \omega, \pi_{\Omega} \}\},Y\}
= \pi^{\sharp}_{\Omega} (i_{X\wedge Y} \omega) \ ,
$$
for all $X$ and $Y \in {\mathcal{X}}(T^*M)$.

$\bullet$ In the notations of Section \ref{MA2} we consider the function $S \in {\mathcal F}(T^*M)$,
$$
S = \mu + \omega \ ,
$$
where $\omega$ is a closed effective $3$-form on $T^*M$.
Then $\{S,S\} = 0$ and $S$ defines a Courant algebroid structure on
$T(T^*M)
\oplus T^*(T^*M)$.

$\bullet$ We consider the Hitchin endomorphism
$H_{\omega}$ and the Poisson bivector $\pi_{\Omega}$ inverse of
$\Omega$. Then
$ H_{\omega}\circ\pi^{\sharp}_{\Omega}
= \pi^{\sharp}_{\Omega}\circ H_{\omega}^*$.
We obtain the following result, the second part of which can be viewed as a
corollary of Theorem 2.5 of \cite{A}.

\begin{theorem}\label{theoremMA3}
Let $\omega$ be the $3$-form on $T^*M$ defined by Formula
\eqref{hess} or \eqref{SL} or \eqref{psSL}, and let
${\mathbb J}_{\omega}$ be the endomorphism of
$T(T^*M) \oplus T^*(T^*M)$ defined by
$$
    {\mathbb J}_{\omega} =
    \begin{pmatrix}
    H_{\omega}& \pi^{\sharp}_{\Omega}\\
    0& -H_{\omega}^*\\
    \end{pmatrix} \ ,
$$
where $H_{\omega}$ is the
Hitchin endomorphism.
If $\lambda_\omega=-1$ (resp., $\lambda_\omega= 1$),
the endomorphism ${\mathbb J}_{\omega}$ is a generalized complex
structure (resp., generalized product structure) on $(T(T^*M) \oplus
T^*(T^*M), \mu + \omega)$.
The triple $(\pi_{\Omega},H_{\omega}, \omega)$ is a Poisson
Nijenhuis structure with background on the manifold $T^*M$.
\end{theorem}

\noindent{\it Proof} 
When $H_{\omega}^2 = - {\rm Id}$ (resp., $+ {\rm{Id}}$), the endomorphism
  ${\mathbb J}_{\omega}$
is a generalized almost complex (resp.,
generalized almost product) structure. Because $\omega$ has constant
coefficients, these structures are integrable \cite{BTh}.
Expressing the vanishing of the torsion of 
${\mathbb J}_{\omega}$ and expanding the terms of $\{\{\pi_\Omega+
H_\omega,\mu+ \omega \},\pi_\Omega+ H_\omega\} + \{H_\omega^2, \mu +
\omega\}$, 
we obtain, as in \cite{A}, 
$$
{\mathrm{ad}}_{\pi_\Omega}^2(\mu)=0,\quad
\{{\mathrm{ad}}_{\pi_\Omega}(\mu),H_\omega\} -
{\mathrm{ad}}_{\pi_\Omega}(\{H_\omega,\mu\})
={\mathrm{ad}}_{\pi_\Omega}^2(\omega) \ ,
$$
$$
\{\{H_{\omega},\mu \},H_{\omega}\} +
\{{\mathrm{ad}}_{\pi_\Omega}(\omega),H_{\omega}\}  -
{\mathrm{ad}}_{\pi_\Omega}(\{H_{\omega},\omega\}) + \{H_{\omega}^2,\mu\} = 0 \ ,
$$
$$
\{\{H_{\omega},\omega\},H_{\omega}\} + \{H^2_\omega, \omega\} = 0 \ .
$$
Replacing $H_{\omega}^2$ par  ${\rm Id}$ or $- {\rm Id}$, we find that 
the quadruple $(\pi_\Omega, H_\omega, 0, \omega)$ is a Poisson quasi-Nijenhuis 
structure with background, {\it{i.e.}}, the triple $(\pi_\Omega, H_\omega , \omega)$ is
a Poisson-Nijenhuis structure with background in the sense of \cite{A}. 
\hfill $\square$

\subsection{Generalized Calabi-Yau structures}

Theorem \ref{theoremMA3} answers a natural question: what is the
relation between the generalized Calabi-Yau structures in the sense of
Hitchin \cite{H} or Gualtieri \cite{Gua} and the generalized Calabi-Yau
structures introduced by Banos \cite{BTh}
in his study of 
Monge-Amp\`ere structures?

The generalized Calabi-Yau structures in the sense of Hitchin are
special generalized complex structures. According to the definition of
M. Gualtieri \cite{Gua} (which is slightly different
from Hitchin's \cite{H}), a
generalized Calabi-Yau manifold is a
manifold with a generalized complex structure and trivial canonical class.
Theorem \ref{theorem8} shows that the Calabi-Yau
Monge-Amp\`ere structures in the sense of Banos are
generalized c.p.s. (complex, product or subtangent) structures in the sense of
Grabowski \cite{G} and Vaisman \cite{V3}.
Equations
\eqref{hess}, \eqref{SL},
\eqref{psSL} define generalized Calabi-Yau structures on $T^*M$
in the sense of Banos.
The Monge-Amp\`ere structures \eqref{SL} and \eqref{psSL}
(called {\it{special Lagrangian}} and {\it{pseudo-special Lagrangian}},
respectively) define
generalized Calabi-Yau structures
in the sense of Gualtieri, while that of \eqref{hess},
where $H_{\omega}^2= {\mathrm{Id}}$, does not since it
corresponds to a generalized product
structure. In the case of \eqref{SL}, we obtain
the canonical Calabi-Yau 
structure on $T^*({\mathbb R}^3)= {\mathbb C}^3$ with the complex
structure $H_{\omega}$, satisfying $H_{\omega}^2= - {\mathrm{Id}}_{{\mathbb R}^6}$.

\bigskip

\noindent{\it Acknowledgments}
V.R. acknowledges the following grants which partially
supported his research during the composition of this paper:
ANR GIMP 2005-2008, RFBR 08-01-00667, and the
INFN-RFBR project ``Einstein''. His special thanks
go to SISSA for the warm hospitality it extended to him and
to U. Bruzzo for numerous discussions.

The authors thank the referee for his comments.

\bigskip

Yvette Kosmann-Schwarzbach, 

Centre de Math\'ematiques Laurent Schwartz, \'Ecole Polytechnique, 

91128 Palaiseau, France

E.mail: yks@math.polytechnique.fr

\medskip
      
Vladimir Rubtsov, 

D\'epartement de Math\'ematiques, Universit\'e d'Angers, 2 Boulevard Lavoisier, 

49045 Angers, France

E.mail: Volodya.Roubtsov@univ-angers.fr


\begin{thebibliography}{40}

\bibitem{A}
P. Antunes, Poisson quasi-Nijenhuis structures with background,
{\it{Lett. Math. Phys.}} {\bf 86} (2008), no. 1, 33-45.

\bibitem{B}
B. Banos,
Monge-Amp\`ere equations and generalized complex geometry -- the
two-dimensional case, {\it J. Geom. Phys.} {\bf 57} (2007), no. 3, 841-853.

\bibitem{BTh}
---, Op\'erateurs de Monge-Amp\`ere symplectiques en dimensions 3 et 4, Th\`ese de Doctorat,
Universit\'e d'Angers, 2002.

\bibitem{CGM1}
J. F. Cari\~nena, J. Grabowski, G.  Marmo,
Contractions:
Nijenhuis and Saletan tensors for general algebraic structures,
{\it J. Phys.} A  {\bf 34}  (2001),  no. 18, 3769-3789.


\bibitem{CGM2}
---,
Courant algebroid and Lie bialgebroid contractions, {\it J. Phys.} A
{\bf 37} (2004),  no. 19, 5189-5202.



\bibitem{C}
M. Crainic, Generalized complex structures and Lie brackets,
arXiv:math/0412097v2 [math.DG].

\bibitem{DF}
P. A. Damianou and R. L. Fernandes, 
Integrable hierarchies and the modular class, {\it Ann. Inst. Fourier} {\bf{58}} (2008), no. 1, 
107-139.
 
\bibitem{D}
I. Ya. Dorfman, {\it{Dirac Structures and Integrability of Nonlinear
    Evolution Equations}}, John Wiley, 1993.

\bibitem{FF}
B. Fuchssteiner, A. S. Fokas,
Symplectic structures, their B\"{a}cklund transformations
and hereditary symmetries, {\it Phys.} D  {\bf 4}  (1981/82), no. 1, 47-66.


\bibitem {G}
J. Grabowski,
Courant-Nijenhuis tensors and generalized geometries, in
Groups, geometry and physics, {\it Monogr. Real
Acad. Cienc. Exact. Fis.-Qu\'{\i}m. Nat. Zaragoza} {\bf 29}
(2006), 101-112.


\bibitem {GU}
J. Grabowski and P. Urbanski,
Lie algebroids and Poisson-Nijenhuis structures, in Quantization,
deformations and coherent states (Bia\l owie\.za, 1996),
{\it {Rep. Math. Phys.}}  {\bf 40}  (1997), 195-208.

\bibitem{Gua}
M. Gualtieri, Generalized complex geometry, arXiv:math/0703298.

\bibitem {H0}
N. Hitchin, The geometry of three-forms in six dimensions,
\emph{J. Differential Geometry} {\bf 55} (2000), 547-576.

\bibitem{H}
---, Generalized Calabi-Yau
  manifolds, \emph{Q. J. Math.} {\bf 54} (2003), no. 3, 281-308. 

\bibitem{yks92}
Y. Kosmann-Schwarzbach,
  Jacobian  quasi-bialgebras  and  quasi-Poisson  Lie  groups,
in \emph{Mathematical Aspects of Classical Field Theory}, M. Gotay,
J. E. Marsden and V. Moncrief, eds.,
Contemp. Math. 132, 1992, 459-489.

\bibitem{yks1996}
---,
From Poisson algebras to Gerstenhaber algebras,
{\it{Ann. Inst. Fourier}} {\bf 46} (1996), 1243-1274.

\bibitem{yks96}
---, The Lie bialgebroid of a Poisson-Nijenhuis
manifold, {\it{  Lett. Math. Phys.}} {\bf 38}  (1996),  421-428.

\bibitem{yks2004}
---,
Derived brackets,  {\it{  Lett. Math. Phys.}}  {\bf 69}  (2004),
61-87.

\bibitem{PM232}
---,
Quasi, twisted, and all that... in Poisson geometry and Lie
algebroid theory, in {\it   The Breadth of Symplectic and Poisson Geometry},
J.~E.~Marsden and T. Ratiu, eds.,
Progr. Math. 232, Birkh\"auser, 2005,
363-389.

\bibitem{yks2007}
---, Poisson and symplectic functions
in Lie algebroid theory, to appear in the Festschrift for Murray
Gerstenhaber and Jim Stasheff, A. Cattaneo, A. Giaquinto and P. Xu,
eds., Prog. Math., Birkh\"auser, 2009, arXiv:0711.2043.

\bibitem{KM}
Y. Kosmann-Schwarzbach and F. Magri,
Poisson-Nijenhuis  structures,  {\it {Ann.  Inst.  Henri
Poincar{\'e}}} A {\bf 53} (1990), 35-81.

\bibitem{KM1}
---,
On the modular classes of Poisson-Nijenhuis manifolds,
preprint, math.SG/0611202.  

\bibitem{KS}
B. Kostant and S. Sternberg,
Symplectic reduction, BRS cohomology, and infinite-dimensional 
Clifford algebras,  
{\it {Ann. Physics}} {\bf 176}  (1987),  no. 1, 49--113.

\bibitem{K}
I. S. Krasil'shchik, Schouten brackets and canonical algebras, in {\it Global
Analysis -- Studies and Applications III}, 
Lecture Notes Math. 1334, Springer,1988, 79-110.

\bibitem{KLR}
A. Kushner, V. Lychagin, V. Rubtsov,
{\it Contact Geometry and Non-Linear Differential Equations}, Encyclopedia
of Mathematics and its Applications 101, Cambridge University Press,
 2007.

\bibitem{LR}
P. Lecomte and C. Roger,
Modules et cohomologies des big\`ebres de Lie,  
{\it {C. R. Acad. Sci. Paris S\'er. I Math.}} {\bf 310}  (1990),
no. 6, 405--410. 

\bibitem{LM}
P. Libermann and C.-M. Marle, {\it Symplectic Geometry and Analytical
Mechanics}, Reidel, 1987.


\bibitem{L}
U. Lindstr\"om, R. Minasian, A. Tomasiello, M. Zabzine,
Generalized complex manifolds and supersymmetry,
\emph{Commun. Math. Phys.} {\bf 257} (2005), 235-256.

\bibitem{Ly}
V. V. Lychagin,
Contact geometry and second-order nonlinear differential equations,
\emph{ Uspekhi Mat. Nauk} {\bf 34}  (1979), no. 1(205), 137-165. English
translation,
\emph{Russian Math. Surveys} {\bf 34} (1979), no. 1, 149-180.

\bibitem{LyR}
V. V. Lychagin, V. N. Rubtsov, I. V. Chekalov, A classification of
Monge-Amp\`ere equations,
\emph{Ann. Sci. \'Ecole Norm. Sup.} 4$^{\rm e}$ s\'erie, {\bf 26} (1993),
no. 3, 281-308.


\bibitem{MM}
F. Magri and C. Morosi, A geometrical characterization of integrable
hamiltonian systems through the theory of Poisson-Nijenhuis
manifolds, Quaderno S/19, Milan, 1984. Re-issued: Universit\`a di Milano
Bicocca, Quaderno 3, 2008. 
http://home.matapp.unimib.it (Quaderni di Dipartimento/2008-3).

\bibitem{MMR}
F. Magri, C. Morosi, O. Ragnisco, Reduction techniques for
infinite-dimensional Hamiltonian systems: some ideas and applications,
{\it Comm. Math. Phys.} {\bf 99} (1985),  115-140.


\bibitem{R}
D. Roytenberg, Quasi-Lie bialgebroids and twisted Poisson manifolds,
{\it Lett. Math. Phys.} {\bf 61} (2002), 123-137.


\bibitem{R2}
---, On the structure of graded symplectic supermanifolds
and Courant algebroids, in \emph{Quantization, Poisson Brackets and
  Beyond} (Manchester, 2001), T. Voronov, ed., Contemp. Math. 315,
2002, 169-185.


\bibitem{SW}
P. \v Severa and A. Weinstein, Poisson geometry with a 3-form
background, in Noncommutative geometry and string theory
(Yokohama, 2001),
\emph{Progr. Theoret. Phys. Suppl.} {\bf 144}  (2001), 145-154.

\bibitem{SX}
M. Sti\'enon and Ping Xu,
Poisson quasi-Nijenhuis manifolds, {\it Comm. Math. Phys.} {\bf 270} (2007), no. 3,
709-725.

\bibitem{V1}
I. Vaisman, Complementary $2$-forms of Poisson structures,
{\it Compositio Math.} {\bf 101} (1996), 55-75.

\bibitem{V2}
---, A lecture on Poisson-Nijenhuis structures, in \emph{Integrable
Systems and Foliations}, C. Albert, R. Brouzet and J.-P. Dufour, eds.,
  Prog. Math. 145, Birkh\"auser, 1997, 169-185.

\bibitem{V3}
---, Reduction and submanifolds of generalized complex manifolds,
{\it Differential Geom. Appl.} {\bf 25} (2007),  no. 2, 147-166.


\bibitem{V4}
---, private communication (1999).


\bibitem{W}
A. Weil, {\it Introduction \`a l'\'etude des vari\'et\'es
  k\"ahl\'eriennes}, Hermann, 1958.


\end{thebibliography}
\end{document}